\documentclass{amsart}
\usepackage{amscd,amssymb}
\usepackage[dvips,final]{graphics}

\theoremstyle{plain}
\newtheorem{thm}{Theorem}[section]

\newtheorem{lem}[thm]{Lemma}
\newtheorem{cor}[thm]{Corollary}
\theoremstyle{remark}
\newtheorem{rem}[thm]{Remark}
\newcommand{\m}{\phantom{-}}
\newcommand{\x}{\negmedspace}
\newcommand{\N}{\mathbb{N}}
\newcommand{\Z}{\mathbb{Z}}
\newcommand{\R}{\mathbb{R}}
\newcommand{\C}{\mathbb{C}\mkern1mu}
\renewcommand{\H}{\mathbb{H}\mkern1mu}
\newcommand{\Ca}{\mathbb{O}\mkern 1mu}
\newcommand{\RP}{\mathbb{R\mkern1mu P}}
\newcommand{\CP}{\mathbb{C\mkern1mu P}}

\newcommand{\Sph}{\mathbb{S}}

\DeclareMathOperator{\diag}{diag}

\DeclareMathOperator{\id}{id}

\newcommand{\1}{{\mathchoice
{\mathrm 1\mskip-4.2mu\mathrm l}{\mathrm 1\mskip-4.2mu\mathrm l}
{\mathrm 1\mskip-3.9mu\mathrm l}{\mathrm 1\mskip-4.0mu\mathrm l}}}

\DeclareMathOperator{\spann}{span}
\DeclareMathOperator{\dist}{dist}

\newcommand{\Wied}{\Sigma^1}
\newcommand{\smc}{\Sigma^3_0}
\newcommand{\sme}{\Sigma^3_1}
\newcommand{\sma}{\Sigma^3_2}
\newcommand{\bmat}{\left[\begin{smallmatrix}}
\newcommand{\emat}{\end{smallmatrix}\right]}
\newcommand{\bsmat}{\bigl[\begin{smallmatrix}}
\newcommand{\esmat}{\end{smallmatrix}\bigr]}
\newcommand{\Bsmat}{\Bigl[\begin{smallmatrix}}
\newcommand{\Esmat}{\end{smallmatrix}\Bigr]}
\newcommand{\bbsmat}{\biggl[\begin{smallmatrix}}
\newcommand{\eesmat}{\end{smallmatrix}\biggr]}
\newcommand{\BBsmat}{\Biggl[\begin{smallmatrix}}
\newcommand{\EEsmat}{\end{smallmatrix}\Biggr]}
\newcommand{\abs}[1]{\vert #1\vert}

\newcommand{\SO}{\mathrm{SO}}
\newcommand{\OO}{\mathrm{O}}

\newcommand{\SL}{\mathrm{SL}}
\newcommand{\U}{\mathrm{U}}

\newcommand{\Syp}{\mathrm{Sp}}
\newcommand{\Gtwo}{\mathrm{G}_2}

\newcommand{\asyp}{\mathfrak{sp}}

\DeclareMathOperator{\Real}{Re}
\DeclareMathOperator{\Imag}{Im}
\DeclareMathOperator{\Ad}{Ad}
\DeclareMathOperator{\ad}{ad}

\newcommand{\ip}{\langle \,\cdot\,,\,\cdot\,\rangle}
\begin{document}
\title[A minimal Brieskorn 5-sphere in the Gromoll-Meyer sphere]{A minimal Brieskorn
5-sphere in the Gromoll-Meyer sphere and its applications}

\author{Carlos Dur\'an}
\address{IMECC-UNICAMP, Pra\c{c}a Sergio Buarque de Holanda, 651, 
Cidade Universit\'aria - Bar\~ao Geraldo, 
Caixa Postal: 6065 
13083-859 Campinas, SP, Brasil }
\email{cduran@ime.unicamp.br}

\author{Thomas P\smash{\"u}ttmann}
\address{Mathematisches Institut\\Universit\"at Bonn\\
        D-53115 Bonn\\Germany}
\email{puttmann@math.uni-bonn.de}
\thanks{The first author was supported by FAPESP grant 03/016789 and
FAEPEX grant 15406, the second author by a DFG Heisenberg fellowship
and by the DFG priority program SPP~1154 ``Globale Differentialgeometrie''.}

\begin{abstract}
We recognize the Gromoll-Meyer sphere $\Sigma^7$ as the geodesic join of a
simple closed geodesic and a minimal subsphere $\Sigma^5\subset\Sigma^7$,
which can be equivariantly identified with the Brieskorn sphere $W^5_3$.
As applications we in particular determine the full isometry group of $\Sigma^7$,
classify all closed subgroups which act freely, determine the homotopy types of the
corresponding orbit spaces, identify the Hirsch-Milnor involution in dimension $5$
with the Calabi involution of $W^5_3$, and obtain explicit formulas for diffeomorphisms
between the two Brieskorn spheres $W^5_3$ and $W^{13}_3$ and standard
Euclidean spheres.
\end{abstract}


\maketitle

%
%
\section{Introduction}
In 1974 Gromoll and Meyer \cite{meyer} constructed an exotic sphere as biquotient
of the compact group $\Syp(2)$ and thereby the first exotic sphere with nonnegative
sectional curvature. This sphere, $\Sigma^7$, can be regarded as the basic example
of a biquotient in Riemannian geometry and, simultaneously, as the basic example
of an exotic sphere. $\Sigma^7$ is naturally an $\Sph^3$-bundle over $\Sph^4$ and
by choosing two local trivializations of this bundle properly, $\Sigma^7$ is identified
with the Milnor sphere $\Sigma^7_{2,-1}$, which is a generator of the group of
homotopy spheres $\Theta_7 \approx \Z_{28}$ in dimension\,7.
Recently, it was shown that $\Sigma^7$ is actually the only exotic sphere that can
be modeled by a biquotient of a compact Lie group \cite{kapovitch}, \cite{totaro}.

\smallskip

Because of this exceptional status of the Gromoll-Meyer sphere it seems natural
to study the geometry of $\Sigma^7$ in detail. Papers that do this from various
viewpoints are \cite{duran}, \cite{eschenburg}, \cite{paternain}, \cite{yamato},
\cite{wilhelm}, for example.
Here, we investigate $\Sigma^7$ through the interaction between symmetry
arguments, submanifold stratifications, and geodesic constructions.
It is important, however, to note that we do not only consider the Gromoll-Meyer metric
on $\Sigma^7$ but the entire $2$-parameter family of metrics $\ip_{\mu,\nu}$
that are $\OO(2)\times \SO(3)$ invariant by construction.
This family includes the Gromoll-Meyer metric ($\mu = \nu =\tfrac{1}{2}$)
and the pointed wiedersehen metric constructed in \cite{duran}
($\mu = \nu = 1$) but not the metrics of almost positive sectional curvature
obtained in \cite{eschenburg} and \cite{wilhelm}.
Extending the constructions of \cite{duran} and \cite{involutions}
we obtain the following structural information:

\begin{thm}
\label{geodesicjoin}
The Gromoll-Meyer sphere $\Sigma^7$ is the join of a simple closed geodesic
$\Sigma^1$ and a minimal subsphere $\Sigma^5$, which is
$\OO(2)\times \SO(3)$-equivariantly diffeomorphic to the Brieskorn sphere $W^5_3$.
For $\mu=1$ and $\nu > 0$ the distance between $\Sigma^1$ and $\Sigma^5$
is constant $\tfrac{\pi}{2}$ and the join structure is realized by distance minimizing
geodesics from $\Sigma^1$ to $\Sigma^5$.
\end{thm}
This theorem and its applications concern the interplay between the Riemannian geometry of the metrics $\ip_{\mu,\nu}$ on $\Sigma^7$ and the equivariant geometry
of the Brieskorn sphere $W^5_3$.
Recall that the Brieskorn sphere $W^5_3$ is the submanifold of $\C^4$ defined by
\begin{gather*}
  z_0^3 + z_1^2 + z_2^2 + z_3^2 = 0,\\
  \abs{z_0}^2 + \abs{z_1}^2 + \abs{z_2}^2 + \abs{z_3}^3 = 1
\end{gather*}
and that there is a natural action of $\OO(2)\times\SO(3)$ on $W^5_3$
(see section\,\ref{briesid}).
On the one hand, it is perhaps not surprising that $W^5_3$
plays a central role for the geometry of~$\Sigma^7$
if one recalls that $\Sigma^7$ is diffeomorphic to $W^7_{6j-1,3}$ for any
$j \in \{1,9,17,\ldots\}$ (see \cite{brieskorn}).
Here, $W^7_{6j-1,3}\subset \C\oplus \C^4$ is defined by the equations
\begin{gather*}
  u^{6j-1} + z_0^3 + z_1^2 + z_2^2 + z_3^2 = 0,\\
  \abs{u}^2 +  \abs{z_0}^2 + \abs{z_1}^2 + \abs{z_2}^2 + \abs{z_3}^3 = 1.
\end{gather*}
On the other hand, while $\Sigma^5$ and $W^5_3$ are
$\OO(2)\times\SO(3)$-equivariantly diffeomorphic,
the ambient spaces $\Sigma^7$ and $W^7_{6j-1,3}$
are {\em never} even $\SO(3)$-equivariantly diffeomorphic
(see Corollary\,\ref{eqhombries}).

It is important to note that $W^5_3$ is not equivariantly diffeomorphic to
$\Sph^5$ with any linear $\OO(2)\times \SO(3)$-action and this holds true
if one restricts from $\OO(2)\times \SO(3)$ to the subgroup
$\OO(3) = \{\pm 1\}\times \SO(3)$. This follows from the classification
theorems of J\"anich and Hsiang/Hsiang (see \cite{bdbook}, \cite{mayer}).
On the other hand, these theorems imply that there exist
$\SO(3)$-equivariant diffeomorphisms $\Sph^5 \to W^5_3$
where $\SO(3)$ acts diagonally on $\Sph^5\subset\R^3\times\R^3$.

This brings us to the first application of Theorem~\ref{geodesicjoin}.
Using the geodesic join structure we derive an explicit formula
for an $\SO(3)$-equivariant diffeomorphism $\Sph^5 \to W^5_3$.
This non-trivial formula can be verified by a straightforward computation
and immediately carries over to dimension~$13$:

\begin{thm}
\label{bi}
Formula (\ref{identbrieskorn}) in section \ref{briesid} provides an
$\SO(3)$-equivariant diffeomorphism $\Sph^5 \to W^5_3$ and a
$\Gtwo$-equivariant diffeomorphism $\Sph^{13}\to W^{13}_3$.
\end{thm}
According to our knowledge this is the first explicit formula for
diffeomorphisms between standard spheres and Brieskorn spheres
$W^{2n-1}_d$ with $n > 2$ and odd $d > 1$.

Related to Theorem\;\ref{bi} is Theorem\;\ref{nonlinear} where we provide
nonlinear actions of $\OO(2)\times\SO(3)$ and $\OO(2)\times\Gtwo$
on the Euclidean spheres $\Sph^5$ and $\Sph^{13}$ that are
equivalent to the $\OO(2)\times \SO(3)$-action on $W^5_3$ and to the
$\OO(2)\times\Gtwo$-action on $W^{13}_3$, respectively.
Various models existed for these actions previously (see \cite{bredon})
but only on manifolds that were inexplicitly diffeomorphic to $\Sph^5$
and~$\Sph^{13}$.

\smallskip

In order to explain the second application of Theorem\,\ref{geodesicjoin}
we have to digress briefly into the history of exotic involutions of spheres.
A fixed point free involution of the standard sphere is called exotic,
if the quotient is not diffeomorphic to the real projective space.
The first examples of such involutions were given by Hirsch and
Milnor \cite{hirsch}. They considered the exotic Milnor sphere
$\Sph^3\cdots \Sigma^7_{2,-1} \to \Sph^4$ and the involution
of $\Sigma^7_{2,-1}$ induced by the antipodal map on the $\Sph^3$-fibers,
detected invariant subspheres of dimensions~$5$ and~$6$
in $\Sigma^7_{2,-1}$, and proved that the restrictions of the involution
of $\Sigma^7$ yield exotic involutions of these subspheres.
The next example of an exotic involution was given by Calabi (unpublished,
see \cite{bredon}) who showed that the involution
of $W^5_3$ given by the map $(z_0,z')\mapsto (z_0,-z')$ is exotic.
In \cite{mayer}, the Calabi involution was identified with an involution
constructed by Bredon and this in turn was identified with the
Hirsch-Milnor involution by Yang \cite{yang}. The latter identification,
however, turned out to be incorrect (see a footnote in \cite{scharlemann}),
so that an explicit identification between the Hirsch-Milnor involution and
the Calabi involution was still missing.

In \cite{involutions} it was shown that the Hirsch-Milnor involutions
are induced by the action of $-\1 \in\Syp(2)$ on $\Sigma^7 = \Sigma^7_{2,-1}$
and that the invariant subspheres of Hirsch and Milnor are precisely the
sphere $\Sigma^5$ and the intermediate equators
$\Sigma^5\subset \Sigma^6_{\pm A} \subset \Sigma^7$
(see section\,\ref{definition}).
In combination with explicit diffeomorphisms $\Sph^5 \to \Sigma^5$
and $\Sph^6 \to \Sigma^6_{\pm A}$ provided by the geodesic geometry
of $\Sigma^7$ this was used to derive explicit formulas for exotic involutions
of the Euclidean spheres $\Sph^5$, $\Sph^6$, $\Sph^{13}$, and $\Sph^{14}$.
As consequence of Theorem\,\ref{geodesicjoin} we now obtain

\begin{cor}
The equivariant diffeomorphism $\Sigma^5 \to W^5_3$ identifies the
Hirsch-Milnor involution in dimension $5$ with the Calabi involution
of $W^5_3$.
\end{cor}

Since $\Sigma^5/\{\pm \1\}$ and $\Sigma^6_{\pm A}/\{\pm \1\}$ are not
diffeomorphic to real projective spaces it is natural to investigate the metrics
on $\Sigma^5$ and $\Sigma^6_{\pm A}$ induced by the metrics $\ip_{\mu,\nu}$
on $\Sigma^7$. We will show that for none of these metrics $\Sigma^5$ or
$\Sigma^6_{\pm A}$ are totally geodesic in $\Sigma^7$. Moreover,
the induced metrics always have some negative sectional curvatures.

\smallskip

The third application of Theorem\,\ref{geodesicjoin} concerns the
full isometry group of~$\Sigma^7$. As mentioned already in \cite{meyer},
Hsiang showed that the maximum dimension of the isometry group of any
metric on $\Sigma^7$ is $4$ (see \cite{straume} for a proof).
Thus the identity component of the isometry group of $\ip_{\mu,\nu}$ is the group
$\SO(2)\times \SO(3)$. It remains the nontrivial problem to determine the
other components. Of particular interest is to see which finite groups can
act freely on $\Sigma^7$. Recent papers \cite{shankar}, \cite{gsz} dealt
with the analogous problem for the homogeneous and cohomogeneity one manifolds
of positive sectional curvature. Surprisingly, it turned out that sometimes noncyclic
groups can act freely on these spaces.
For the cohomogeneity~$3$ metrics $\ip_{\mu,\nu}$ on $\Sigma^7$ the
structural information of Theorem\,\ref{geodesicjoin} can be used to reduce
the problem to the corresponding problem for the induced cohomogeneity one
metrics on $\Sigma^5$. This latter problem can be solved with the help of some
curvature computations.
\begin{thm}
\label{freeness}
The group $\OO(2) \times \SO(3)$ is the full isometry group of the metrics
$\ip_{\mu,\nu}$ on $\Sigma^7$ and on $\Sigma^5$. Any subgroup that acts
freely on either $\Sigma^7$ or $\Sigma^5$ is a finite cyclic group.
Conversely, for any $m\in \Z$ the group $\Z_m$ acts freely and isometrically on
$\Sigma^7$ and on $\Sigma^5$, even in several non-conjugate ways for a
fixed $m > 2$.
\end{thm}

All the $\Z_m$-quotients of $\Sigma^7$ inherit nonnegative sectional
curvature from the Gromoll-Meyer metric $\ip_{\frac{1}{2},\frac{1}{2}}$.
(It is interesting, however, to note that for $m > 2$ the known metrics
with almost positive sectional curvature on $\Sigma^7$ \cite{eschenburg},
\cite{wilhelm} are not invariant with respect to the $\Z_m$-actions.)
In the case of $\Sigma^5$ none of the metrics
$\ip_{\mu,\nu}$ has nonnegative curvature but by the Grove-Ziller construction
\cite{grove} there exist $\OO(2)\times\SO(3)$-invariant metrics
on $\Sigma^5$ with $K \ge 0$.

\begin{cor}
\label{exoticlens}
For any $m$ not divisible by $6$ and any two integers $p$, $q$ with
$p\neq 0$, $3p-q\neq 0$, $3p+q \neq 0$ such that $m$ is relatively
prime to $p$, $3p-q$, and $3p+q$ there is a $7$-manifold with
$K \ge 0$ that is homotopy equivalent to the lens space $L^7_m(p,p,3p-q,3p+q)$
but not diffeomorphic to any standard lens space.
If $m$ is even there also exists a $5$-manifold with $K \ge 0$ that is homotopy
equivalent to the lens space $L^5_m(p,3p-q,3p+q)$ but not diffeomorphic to
any standard lens space.
\end{cor}

In dimension~5 the case $m=2$ was already covered in \cite{grove}.
Apart from this these seem to be the first known exotic homotopy lens spaces
with $K \ge 0$. Exotic lens spaces with positive Ricci and almost nonnegative
sectional curvature in higher dimensions were found by Schwachh\"ofer and
Tuschmann \cite{schwachtusch}.

The non-trivial part of Corollary\,\ref{exoticlens} is to determine the homotopy
type of the free $\Z_m$-quotients of $\Sigma^7$ and $\Sigma^5$.
It is well-known that the orbit spaces of free $\Z_m$-actions on homotopy
spheres are homotopy equivalent to lens spaces (see \cite{browder}).
For a concretely given action on a homotopy sphere, however, there is no
canonical tractable way to determine the homotopy type of the quotient.
In our case we will follow an idea of Orlik \cite{orlik} and construct
branched coverings $\Sigma^5\to \Sph^5$ that can be extended
by the join structure of Theorem\,\ref{geodesicjoin} to (continuous)
branched coverings $\Sigma^7 \to \Sph^7$.
The essential property of these branched coverings is that they are
$\OO(2)\times \SO(3)$-equivariant where $\OO(2)\times \SO(3)$ acts
linearly on $\Sph^7$ and $\Sph^5$. The target spaces of the induced maps
$\Sigma^7/\Z_m \to \Sph^7/\Z_m$ and $\Sigma^5/\Z_m\to \Sph^5/\Z_m$
are thus standard lens spaces and this allows us to determine the
homotopy type of $\Sigma^7/\Z_m$ and $\Sigma^5/\Z_m$.

\smallskip

The fourth application of Theorem\,\ref{geodesicjoin} resides in
determining the structure of fixed point sets of isometries of $\Sigma^7$.
Fixed point sets of isometries are useful to understand the geometry
of Riemannian manifolds since each connected component is a
totally geodesic submanifold. In particular, they provide significant
curvature information since the extrinsic and intrinsic sectional curvature
of a plane tangent to a totally geodesic submanifold are equal.
In a general biquotient it is fairly difficult if not impossible to determine
the structure of all fixed point sets. In $\Sigma^7$, however,
we can employ Theorem\,\ref{geodesicjoin} to determine the metric structure
of all fixed point sets in a very geometric way (see section\,\ref{totgeod}).
It turns out that all fixed point sets with dimension $> 1$ are congruent to
one of three $3$-spheres $\smc$, $\sme$, $\sma$, to a real projective space $P^3$,
or a lens space $L^3\approx \Sph^3/\Z_3$.
It is interesting to see how the biquotient structure of $\Sigma^7$ causes
$\smc$ and $L^3$ to have more intrinsic than extrinsic isometries:
Both are intrinsically homogeneous although they only inherit a cohomogeneity
one action from the $\OO(2)\times \SO(3)$ action on $\Sigma^7$.
The fact that $L^3$ and $P^3$ are fixed point sets with non-trivial fundamental group
shows how much the geometry of $\Sigma^7$ differs from the geometry
of the standard sphere $\Sph^7$ with constant curvature.
The induced metrics $\ip_{\mu,\nu}$ on $\sma$ are properly of cohomogeneity~$2$.
They are remarkable since the curvature tensor looks like the metrics
would be of cohomogeneity~$1$ and there is no obvious deformation to the
constant curvature metric through metrics with this property.

\smallskip

The authors would like to thank Uwe Abresch and A.~Rigas for many
useful discussions. We also thank Luigi Verdiani whose Maple applets
allowed to cross-check some of the curvature computations in a very
efficient way.

\bigskip

\section{The geodesic join description of the Gromoll-Meyer sphere}
\label{definition}
In this section we review and extend some of the constructions of \cite{duran}
and \cite{involutions}. In particular, we use a $1$-parameter family of
metrics on the Gromoll-Meyer sphere $\Sigma^7$ with the pointed wiedersehen
property along a circle $\Sigma^1$ to recognize $\Sigma^7$ as the geodesic join
of $\Sigma^1$ and a minimal subsphere $\Sigma^5\subset \Sigma^7$.

Let $\Sph^3$ denote the unit sphere in the quaternions $\H$ and $\Syp(2)$ the
group of $2\times 2$ quaternionic matrices $A$ such that $\bar A^t\cdot A = \1$.
On $\Syp(2)$ we consider the class of left invariant and $\Syp(1)\times \Syp(1)$
right invariant Riemannian metrics. After rescaling, these metrics correspond
precisely to the $\Ad_{\Syp(1)\times\Syp(1)}$-invariant inner products
\begin{gather*}
  \langle \bsmat x_1 & -\bar y_1\\ y_1 & \m z_1\esmat,
    \bsmat x_2 & -\bar y_2\\ y_2 & \m z_2\esmat\rangle_{\mu,\nu}
  = \Real \bigl( \mu \, \bar x_1 x_2 + \bar y_1 y_2 + \nu \,\bar z_1 z_2 \bigr)
\end{gather*}
on the Lie algebra $\asyp(2)$. The standard biinvariant metric on $\Syp(2)$
is $\ip_{\frac{1}{2},\frac{1}{2}}$. This metric has nonnegative sectional curvature and
it follows from Cheeger's construction \cite{cheeger} that all metrics $\ip_{\mu,\nu}$
with $\mu,\nu \le \tfrac{1}{2}$ have nonnegative sectional curvature as well.

Two free isometric actions of $\Sph^3$ on $\Syp(2)$ play a central role in
the rest of the paper: the standard action
\begin{gather*}
  \Sph^3 \times \Syp(2) \to \Syp(2),\quad
  q \bullet A = A \cdot \bsmat 1& 0\\ 0 & \bar q\esmat
\end{gather*}
and the Gromoll-Meyer action \cite{meyer}
\begin{gather*}
  \Sph^3 \times \Syp(2) \to \Syp(2),\quad
  q \star A = q \cdot A \cdot \bsmat \bar q & 0\\ 0 & 1\esmat.
\end{gather*}
Both these actions foliate $\Syp(2)$ by $\Sph^3$-orbits in two different ways:
The orbit space of the standard action $\bullet$ can be naturally identified
with $\Sph^7\subset \H^2$ by restricting a matrix in $\Syp(2)$ to its first column.
The orbit space $\Sigma^7$ of the Gromoll-Meyer action~$\star$
is diffeomorphic to the exotic Milnor sphere $\Sigma^7_{2,-1}$.
The corresponding projection maps are denoted by $\pi_{\Sph^7}: \Syp(2)\to\Sph^7$
and $\pi_{\Sigma^7}: \Syp(2)\to\Sigma^7$.
Throughout this paper both orbit spaces, $\Sph^7$ and $\Sigma^7$, are supposed
to carry metrics induced from $\ip_{\mu,\nu}$ on $\Syp(2)$ by Riemannian submersion.
The metrics on $\Sph^7$ and $\Sigma^7$ will also be denoted by $\ip_{\mu,\nu}$.
Since Riemannian submersions are curvature nondecreasing it is clear that
the sectional curvature of $(\Sigma^7,\ip_{\mu,\nu})$ is nonnegative for
$\mu,\nu \le \tfrac{1}{2}$.

The starting point for our subsequent geometric constructions and considerations
is the following elementary fact: The $\bullet$-orbit and the $\star$-orbit
through any real matrix $A \in \OO(2)\subset \Syp(2)$ are equal as sets
since $A$ commutes with all $q\in\Sph^3$. A~geodesic in $\Syp(2)$ that passes
perpendicularly through the common orbit
\begin{gather*}
  \Sph^3\bullet A = \Sph^3\star A
  = \bigl\{ A \cdot \bsmat 1& 0\\ 0 & \bar q\esmat \;\big\vert\; q \in \Sph^3\bigr\}
\end{gather*}
is perpen\-dicular to all $\bullet$-orbits and all $\star$-orbits and hence
projects to geodesics in both orbit spaces, $\Sph^7$ and $\Sigma^7$.
(Recall that the inner product between the velocity vector field and a Killing field
is constant along a geodesic).
Now, fixing a matrix $A' \in \Sph^3\bullet A = \Sph^3\star A$, we get an identification
between the geodesics in $\Sph^7$ that start at the point $\Sph^3\bullet A$
and the geodesics in $\Sigma^7$ that start at the point $\Sph^3\star A$:
Two such geodesics $\gamma_{\Sph^7}$ and $\gamma_{\Sigma^7}$
correspond to each other if and only if there is a common horizontal lift through~$A'$,
i.e., a geodesic $\tilde \gamma$ in $\Syp(2)$ which starts at $A'$ perpendicularly to
$\Sph^3\bullet A = \Sph^3\star A$ such that
$\gamma_{\Sph^7} = \pi_{\Sph^7}\circ \tilde\gamma$ and
$\gamma_{\Sigma^7} = \pi_{\Sigma^7}\circ \tilde\gamma$.
Since $\Sph^3\bullet A = \Sph^3\star A$ are only equal as sets, this identification
depends on the choice of~$A'$. There is, however, a canonical choice for $A'$
since $\Sph^3\bullet A = \Sph^3\star A$ intersects $\OO(2)$
precisely in the set $\bigl\{A,A\cdot \bmat 1 & \m 0\\ 0 & -1\emat\bigr\}$
and hence contains a unique element $A' \in \SO(2)$.

\begin{figure}[ht]
\begin{center}
\mbox{\scalebox{0.7}{\includegraphics{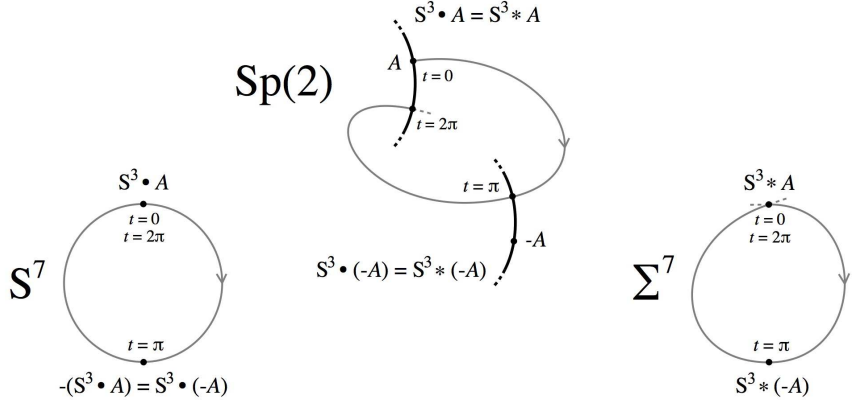}}}
\caption{For $A\in \SO(2)\subset \Syp(2)$ a geodesic in $\Sph^7$ through
$\Sph^3\bullet A$ corresponds precisely to one geodesic in $\Sigma^7$ through
$\Sph^3\star A$ via a common horizontal lift through $A$.}
\end{center}
\end{figure}

This correspondence has an immediate application in the case $\mu=1$
where each left invariant metric $\ip_{1,\nu}$ on $\Syp(2)$ induces the
standard metric on $\Sph^7$: 
All unit speed geodesics of $\Sph^7$ pass through their antipode after
time $\pi$ and return to their starting point after time $2\pi$. This holds
in particular for the geodesics that start at a point $\Sph^3\bullet A$
with $A \in \OO(2)$. Since the antipode of $\Sph^3\bullet A$ is the orbit
$\Sph^3 \bullet (-A) = \Sph^3\star (-A)$, the geodesic correspondence
above implies the following recurrency behavior:

\begin{thm}[see \cite{duran} in the case $\nu = 1$]
\label{wiedersehen}
The unit speed geodesics of $(\Sigma^7,\ip_{1,\nu})$ that start at a point
$\Sph^3\star A$ with $A\in \OO(2)$ all pass through $\Sph^3\star (-A)$
after time $\pi$ and return to $\Sph^3\star A$ after time $2\pi$ (but do not
close smoothly in general).
These geodesics are length minimizing until time $\pi$.
\end{thm}

In accordance with the literature (see e.g. \cite{besse}) the points of
the circle
\begin{gather*}
  \Wied := \{ \Sph^3 \star A \;\vert\; A \in \OO(2) \} \subset \Sigma^7
\end{gather*}
will be called {\em wiedersehen} points. The wiedersehen property allows us
to define natural subspheres of $\Sigma^7$: For $A \in \OO(2)$ the bisector
\begin{gather*}
  \Sigma^6_{\pm A} := \bigl\{ x \in \Sigma^7 \,\big\vert\,
    \dist(x,\Sph^3\star A) = \dist(x,\Sph^3\star(-A)) = \tfrac{\pi}{2}\bigr\}
\end{gather*}
is given by the midpoints of the geodesics that start at $\Sph^3\star A$ and end
at $\Sph^3\star (-A)$. The intersection of all the bisectors $\Sigma^6_{\pm A}$
in $\Sigma^7$ is the set
\begin{gather*}
  \Sigma^5 := \bigcap_{A \in \OO(2)} \Sigma^6_{\pm A}
    = \bigl\{ x \in \Sigma^7 \,\big\vert\,
    \dist\bigl(x,\Wied\bigr) = \tfrac{\pi}{2} \bigr\}.
\end{gather*}

Recall that the join $X\ast Y$ of two spaces $X$ and $Y$ is the quotient
of $X\times Y\times [0,1] / \sim$ where $(x,y,0) \sim (x,y',0)$ and
$(x,y,1) \sim (x',y,1)$ for all $x\in X$ and all $y\in Y$.
For our purposes it is convenient to substitute $[0,1]$ by $[0,\tfrac{\pi}{2}]$.
\begin{cor}
\label{joingeod}
For $\mu = 1$ the Gromoll-Meyer sphere $\Sigma^7$ is the geodesic
join of the circle $\Sigma^1$ and the subsphere $\Sigma^5$ which
have constant distance $\tfrac{\pi}{2}$, i.e., the map
$\Sigma^1 \ast \Sigma^5 \to \Sigma^7$
that maps $(x,y,t)$ to $\gamma(t)$, where $\gamma: [0,\tfrac{\pi}{2}] \to \Sigma^7$
is the unique unit speed geodesic segment from $x$ to $y$, is a homeomorphism.
\end{cor}

\smallskip

The identification of geodesics in $\Sph^7$ that start at $\Sph^3\bullet A$ with the
geodesics of $\Sigma^7$ that start at a point $\Sph^3\star A$ provides an
$\SO(3)$-equivariant homeomorphism between $\Sph^7$ and $\Sigma^7$ that
restricts to a diffeomorphism between $\Sph^7 \smallsetminus (\Sph^3\bullet (-A))$
and $\Sigma^7 \smallsetminus (\Sph^3\star (-A))$.

This diffeomorphism further restricts to diffeomorphisms $S^6_{\pm A} \to \Sigma^6_{\pm A}$
and $S^5 \to \Sigma^5$ where
\begin{align*}
  S^6_{\pm A} &= \bigl\{ \bmat w_1\\ w_2\emat \in \Sph^7\subset \H^2
      \,\big\vert\, \dist\bigl( \bmat w_1\\ w_2\emat, \bmat a_{11}\\ a_{21}\emat \bigr)
        = \dist \bigl( \bmat w_1\\ w_2\emat, -\bmat a_{11}\\ a_{21}\emat \bigr)
        = \tfrac{\pi}{2}\bigr\} \\
     &= \bigl\{ \bmat w_1\\ w_2\emat \in \Sph^7\subset\H^2
      \,\big\vert\, \Real( a_{11} w_1 + a_{21} w_2) = 0 \bigr\}\\
     &= \{ A \cdot \bmat p\\ w\emat \in \Sph^7 \subset \H^2\;\vert\; p \in \Imag \H,\; w \in \H \}
\end{align*}
and
\begin{gather*}
  S^5 = \bigl\{ \bmat p_1\\ \smash[b]{p_2}\emat
    \,\big\vert\, p_1,\smash[b]{p_2} \in \Imag \H,\;
    \abs{p_1}^2 + \abs{\smash[b]{p_2}}^2 = 1\bigr\}.
\end{gather*}
Note that
\begin{align}
\begin{split}
\label{preimages}
  \pi_{\Sigma^7}^{-1}(\Sigma^6_{\pm A}) &= \pi_{\Sph^7}^{-1}(S^6_{\pm A})
  = \{ A \cdot \bmat p & \ast\\ w & \ast \emat \in \Syp(2) \;\vert\; p \in \Imag \H,\; w \in \H \}\\
  \pi_{\Sigma^7}^{-1}(\Sigma^5)\; &= \;\pi_{\Sph^7}^{-1}(S^5) \;\;
  = \{\bmat p_1 & \ast\\ \smash[b]{p_2} & \ast \emat \in \Syp(2) \;\vert\; p_1, \smash[b]{p_2} \in \Imag\H\}
\end{split}
\end{align}
since the two sets on the right hand side are invariant under the $\star$-action.

There are explicit formulas for the horizontal lifts of the relevant geodesics in $\Sph^7$
and hence for the diffeomorphisms $S^6_{\pm A} \to \Sigma^6_{\pm A}$
and $S^5 \to \Sigma^5$: Consider the geodesic
\begin{gather*}
 \gamma_{\bmat p\\ w\emat}(t) =
   \cos t \, \bsmat 1 \\ 0\esmat + \sin t \, \bmat p\\ w\emat
\end{gather*}
in $\Sph^7 \subset \H^2$ that emanates from the north pole with initial velocity
$\bmat p\\ w \emat\in \Sph^6\subset \Imag\H \times \H$.
The unique horizontal lift $\tilde\gamma_{\bmat p\\ w\emat}$ of
$\gamma_{\bmat p\\ w\emat}$ to $\Syp(2)$ with
$\tilde \gamma_{\bmat p\\ w\emat}(0) = \1$ is given by
\begin{gather}
\label{lift}
  \tilde \gamma_{\bmat p\\ w\emat} (t) =
    \cos t \, \bbsmat 1 & 0\phantom{j}\\ 0 & \, 
    \smash[b]{\tfrac{w}{\abs{w}} e^{tp} \tfrac{\bar w}{\abs{w}}} \eesmat
    + \sin t \, \bbsmat p & \; -\smash{e^{tp} \bar w}\\ w & \;
    - \smash[b]{\tfrac{w}{\abs{w}} p e^{tp} \tfrac{\bar w}{\abs{w}}} \eesmat,
\end{gather}
where $e^{p} = \cos \abs{p} + \tfrac{p}{\abs{p}} \sin \abs{p}$
denotes the exponential map of $\Sph^3\subset\H$ at $1$.
Note that for $w = 0$ equation (\ref{lift}) simply becomes
$\tilde\gamma_{\bmat p\\ 0\emat}(t) = \bmat e^{tp} & 0\\ 0 & 1\emat$.
Now the curve $\pi_{\Sigma^7}\circ \tilde\gamma_{\bmat p\\ w\emat}$
is a geodesic of $\Sigma^7$ for all metrics $\ip_{1,\nu}$
and
\begin{gather}
\label{expldiffeo}
  \Sph^6 \to \Sigma^6_{\pm \1}, \quad
  \bmat p\\ w\emat \mapsto
    \pi_{\Sigma^7}\circ \tilde\gamma_{\bmat p\\ w\emat}( \tfrac{\pi}{2} )
\end{gather}
is an analytic diffeomorphism. This diffeomorphism restricts to an
analytic diffeomorphism $\Sph^5 \to \Sigma^5$ for $\Real w = 0$.

\smallskip

In \cite{involutions} it was shown that $\Sigma^5/\{\pm \1\}$ and
$\Sigma^6_{\pm A}/\{\pm \1\}$ are homotopy equivalent but not
diffeomorphic to $\RP^5$ and $\RP^6$, respectively.
We conclude this section with the following observation:
\begin{lem}
\label{geodesic}
Let $A_0 := \bsmat 0 & -1\\ 1 & \m 0\esmat$.
For any $A \in \OO(2)$ the bisector $\Sigma^6_{\pm A\cdot A_0}$ in
$(\Sigma^7, \ip_{1,\nu})$
is geodesic at the two points $\Sph^3\star A$ and $\Sph^3\star (-A)$,
i.e., any geodesic of $(\Sigma^7, \ip_{1,\nu})$ that starts at one of these points
tangentially to $\Sigma^6_{\pm A\cdot A_0}$ is completely contained
in $\Sigma^6_{\pm A\cdot A_0}$.
\end{lem}

\begin{proof}
It suffices to consider the case $A = \1$. By (\ref{preimages}) we have
\begin{gather*}
  \Sigma^6_{\pm A_0} = \pi_{\Sigma^7}\bigl(
    \bigl\{ \bsmat \smash[t]{w'} & \ast\\ \smash[b]{p'} & \ast \esmat \in \Syp(2)
      \;\big\vert\; p' \in \Imag \H,\; w' \in \H \bigr\} \bigr)\,.
\end{gather*}
Form (\ref{lift}) it is now evident that all the geodesics
$\pi_{\Sigma^7}\circ \tilde \gamma_{\bmat p\\ w\emat}$
with $\Real w = 0$ are contained in $\Sigma^6_{\pm A_0}$.
\end{proof}

\begin{cor}
The exotic projective space $\Sigma^6_{\pm A}/\{\pm \1\}$ inherits from $\Sigma^7$
a one parameter family of metrics that are Blaschke at one point.
\end{cor}

\bigskip

\section{The isometry group of the Gromoll-Meyer sphere}
\label{isomgroup}
The $\bullet$-action of $\Sph^3$ on $\Syp(2)$ of the previous section
extends to the action
\begin{gather*}
  \OO(2) \times \Sph^3 \times \Syp(2) \to \Syp(2),\quad
  (A, q) \bullet B =
    A \cdot B \cdot \bsmat 1 & 0\\ 0 & \bar q\esmat.
\end{gather*}
This action is isometric for all metrics $\ip_{\mu,\nu}$ on $\Syp(2)$ and commutes
with the Gromoll-Meyer action $\star$. Hence, it induces an effective isometric action
\begin{gather*}
  \OO(2) \times \SO(3) \times \Sigma^7 \to \Sigma^7,\quad \SO(3) = \Sph^3/\{\pm 1\},
\end{gather*}
on $\Sigma^7$ that will again be denoted by $\bullet$. This action already
appeared in the original paper of Gromoll and Meyer \cite{meyer}. At the end of
this section we will show that $\OO(2)\times \SO(3)$ is the full isometry
group for all the metrics $\ip_{\mu,\nu}$ on $\Sigma^7$. The following simple fact is
fundamental for the rest of the paper. It allows us to investigate geometric properties
of the metrics that $\Sigma^5$ inherits from $\Sigma^7$, it yields an equivariant
diffeomorphism  between $\Sigma^5$ and the Brieskorn sphere $W^5_3$, and
it is the key to determine which isometries act freely on the Gromoll-Meyer sphere.

\begin{lem}
\label{invariant}
The $\bullet$-action of $\OO(2) \times \SO(3)$ on $\Sigma^7$ leaves
$\Sigma^1$ and $\Sigma^5$ invariant.
The induced $\bullet$-action on $\Sigma^5$ is of cohomogeneity one.
\end{lem}

The cohomogeneity one action is studied in detail in section\,\rm\ref{cohomone}.
Note that the $\bullet$-action of $\OO(2)\times \SO(3)$ does not
leave any of the $\Sigma^6_{\pm A}$ invariant. The largest action that preserves
$\Sigma^6_{\pm \1}$ (and also $\Sigma^6_{\pm A_0}$ with
$A_0 = \bmat 0 & -1\\ 1 & \m 0\emat$)
is the restriction of the $\bullet$-action to $\Z_2\times\Z_2 \times \SO(3)$
where $\Z_2\times\Z_2$ is the group of diagonal matrices in the $\OO(2)$-factor.

\begin{cor}
The Gromoll-Meyer sphere $\Sigma^7$ is $\OO(2)\times \SO(3)$-equivariantly
homeomorphic to the join $\Sigma^1 \ast \Sigma^5$.
\end{cor}

\smallskip

In the following three lemmas we use the geodesic constructions from the previous
section to get some immediate structural information about fixed point sets.
Any element  (or any subgroup) of $\OO(2)\times\SO(3)$ either fixes the
entire circle $\Sigma^1$, two antipodal points in $\Sigma^1$, or no points
in $\Sigma^1$ at all.

\begin{lem}
\label{nowiedfix}
If an element $\psi \in \OO(2)\times\SO(3)$ does not have a fixed point
in $\Wied$ then its fixed point set $(\Sigma^7)^{\psi}$ is completely
contained in $\Sigma^5$.
\end{lem}

\begin{proof}
It suffices to consider any of the metrics $\ip_{1,\nu}$. By Lemma\,\ref{invariant} 
the isometry $\psi$ maps $\Wied$ and $\Sigma^5$ to themselves. 
Through any point $p \in \Sigma^7$ outside $\Wied \cup \Sigma^5$ there is a unique
geodesic segment from $\Wied$ to $\Sigma^5$ with length~$\tfrac{\pi}{2}$.
If $p$ is fixed by the isometry $\psi$ this segment is fixed pointwise as well.
\end{proof}

\begin{cor}
\label{minimal}
For all the metrics $\ip_{\mu,\nu}$ the 5-sphere $\Sigma^5$ is a minimal submanifold
of $\Sigma^7$ and of each $\Sigma^6_{\pm A}$.
\end{cor}
\begin{proof}
All principal isotropy groups of the $\bullet$-action on $\Sigma^5\subset\Sigma^7$
are conjugate to the subgroup $H\subset \OO(2)\times\SO(3)$ determined
in Lemma\,\ref{isotropy}. The union $\Sigma^7_{(H)}$ of all orbits of type $(H)$ in
$\Sigma^7$ (i.e., the set of all points whose isotropy group is conjugate to $H$) is
a perhaps disconnected open minimal submanifold of $\Sigma^7$ (see \cite{lawson}).
Any subgroup conjugate to $H$ contains an element of the form $(-\1,\pm q)$.
All elements of this form act on $\Sigma^1$ by the antipodal map. Thus,
Lemma\,\ref{nowiedfix} implies that $\Sigma^7_{(H)}$ is contained in $\Sigma^5$.
Now, $\Sigma^5$ is apparently the closure of $\Sigma^7_{(H)}$ and hence minimal.
An analogous argument using the $\Z_2\times\Z_2\times\SO(3)$ action
shows that $\Sigma^5$ is minimal in~$\Sigma^6_{\pm \1}$.
\end{proof}

\begin{lem}
\label{twofixwied}
If $\psi\in \OO(2)\times\SO(3)$ fixes precisely two antipodal points
$\Sph^3 \star(\pm A) \in \Wied$ with $A\in \OO(2)$ then $(\Sigma^7)^{\psi}$
is contained in
$\Sigma^6_{\pm A\cdot A_0}$ with $A_0 = \bmat 0 & -1\\ 1 & \m 0\emat$.
Moreover, $(\Sigma^7)^{\psi}$ is a suspension of $(\Sigma^5)^{\psi}$
from the two points $\Sph^3 \star(\pm A)$. In particular, $(\Sigma^7)^{\psi}$
and $(\Sigma^5)^{\psi}$ are both diffeomorphic to spheres.
\end{lem}

\begin{proof}
It suffices to consider any of the metrics $\ip_{1,\nu}$. Between any two points
$p,q \in \Sigma^7$ with $p\in \Sigma^1$ and
$q \not \in \Sigma^1$ there exists a unique minimizing geodesic segment
from $p$ to $q$. If $\psi$ fixes $p$ and $q$ then it fixes the geodesic as well.
Note that $(\Sigma^7)^{\psi}$ is either empty or odd dimensional
since $\psi$ is orientation preserving. (As a generator of the group
$\Theta_7\approx\Z_{28}$ of homotopy spheres, $\Sigma^7$ does
not admit orientation reversing diffeomorphisms.)
\end{proof}

\begin{lem}
\label{fullwiedfix}
If $\psi\in \OO(2)\times\SO(3)$ fixes all points in the circle $\Wied$
then $(\Sigma^7)^{\psi}$ is either equal to $\Wied$, or $(\Sigma^7)^{\psi}$
is the join of $\Wied$ and $(\Sigma^5)^{\psi}$ and hence the suspension
of $(\Sigma^6_{\pm A})^{\psi}$ from any two antipodal points
$\Sph^3 \star(\pm A) \in\Wied$.
In particular, $(\Sigma^7)^{\psi}$, $(\Sigma^6_{\pm A})^{\psi}$,
and $(\Sigma^5)^{\psi}$ are diffeomorphic to spheres.
\end{lem}
\begin{proof}
Similar to the proofs of Lemma\,\ref{nowiedfix} and Lemma\,\ref{twofixwied}.
\end{proof}

The three lemmas above show that the topologically interesting fixed point
sets are all contained in $\Sigma^5$. In section\,\ref{totgeod} we will study
the induced metrics on all existing fixed point sets.

In Lemma\,\ref{invariant} it was shown that all elements of $\OO(2)\times\SO(3)$
map $\Wied$ and $\Sigma^5$ to themselves. The same is true for any isometry
of $(\Sigma^7,\ip_{\mu,\nu})$:

\begin{lem}
\label{invariantspheres}
Every isometry of $(\Sigma^7, \ip_{\mu,\nu})$ maps $\Wied$ and $\Sigma^5$
to themselves.
\end{lem}

\begin{proof}
The maximum dimension of any compact differentiable transformation group on
$\Sigma^7$ is $4$ (see \cite{straume}). The $\bullet$-action of
$G = \OO(2)\times \SO(3)$ on $\Sigma^7$ is effective.
Hence, the subgroup $G_0 = \SO(2)\times \SO(3)$ of $G$
is the identity component of the isometry group $\tilde G$ of $\Sigma^7$.
Let $\psi \in \tilde G$ be any isometry of $\Sigma^7$. Since $G_0$ is a normal subgroup
of~$\tilde G$ conjugation by $\psi$ on $\tilde G$ maps $\SO(3)$ to itself.
Hence, $\psi$ maps the fixed point set of $\SO(3)$ in $\Sigma^7$ to itself.
This fixed point set is precisely the circle $\Wied$ of wiedersehen points.
Moreover, $\psi$ also maps $G_0$-orbits diffeomorphically to $G_0$-orbits.
It thus follows from the isotropy groups like in the proof of Corollary\,\ref{minimal}
that $\psi$ maps $\Sigma^5$ to itself.
\end{proof}

\begin{thm}
The isometry group of $(\Sigma^7,\ip_{\mu,\nu})$ is the group
$\OO(2)\times \SO(3)$.
\end{thm}

\begin{proof}
We need some geometric facts from the following sections for the proof.
In Lemma\,\ref{isomsigmafive} it will be shown that $\OO(2)\times\SO(3)$
is the isometry group of $(\Sigma^5,\ip_{\mu,\nu})$.
By the previous lemma, every isometry of $(\Sigma^7, \ip_{\mu,\nu})$ maps
$\Sigma^5$ to itself. It suffices to show that the only isometry of $\Sigma^7$ that
fixes $\Sigma^5$ pointwise is the identity. Let $\psi$ be such an isometry.
The fixed point set of $\psi$ is the disjoint union of totally geodesic submanifolds.
Consider the component $M$ that contains $\Sigma^5$. Clearly,
$\Sigma^2 \subset M\cap \sma$ (see section\,\ref{totgeod}).
Since $\Sigma^2$ is not totally geodesic in $\sma$ it follows that $\sma \subset M$.
Now consider the congruent copy $\tilde \sma$ of $\sma$ given by the fixed
point set of the isometry
$\bigl( \bsmat -1 & \m 0\\ \m 0 & \m 1\esmat, \pm k\bigr)\bigr\}$ in $\Sigma^7$.
By the same argument as above $\tilde \sma$ is also contained in $M$.
The inclusion $\Sigma^5 \cup \sma \cup \tilde \sma \subset M$
implies that $\dim M = 7$ (inspect the tangent spaces along the
normal geodesic) and hence that $\psi = \id_{\Sigma^7}$.
\end{proof}

\bigskip\bigskip\bigskip

\section{The cohomogeneity one action on $\Sigma^5$}
\label{cohomone}
We will now study the cohomogeneity one action $\bullet$ on $\Sigma^5$ in detail.
The essential technical step is to find a normal geodesic, i.e., a geodesic that
crosses all $\bullet$-orbits perpendicularly.
Recall from (\ref{preimages}) that
\begin{gather*}
  \Sigma^5\; = \;\pi_{\Sph^7} \bigl(
  \{ \bmat p_1 & \ast\\ \smash[b]{p_2} & \ast \emat \in \Syp(2)
    \;\vert\; p_1, \smash[b]{p_2} \in \Imag\H\}\bigr).
\end{gather*}
We will show that the curve
$\alpha(s) = \pi_{\Sigma^7}(\tilde \alpha(s))$ with
\begin{gather*}
  \tilde \alpha(s) = \bsmat j\cos s & \m k \sin s\\ k \sin s & \m j\cos s\esmat
\end{gather*}
is such a normal geodesic and compute the isotropy groups along this geodesic
and the induced Riemannian metrics on the principal orbits. Finally, we will
show that $\OO(2)\times\SO(2)$ is the full isometry group of $(\Sigma^5,\ip_{\mu,\nu})$.

\begin{lem}
The curve $\tilde \alpha$ intersects all $\star$-orbits in $\Syp(2)$ perpendicularly,
i.e., $\tilde \alpha$ is horizontal with respect to the submersion
$\pi_{\Sigma^7}: \Syp(2) \to \Sigma^7$.
\end{lem}

\begin{proof}
For all the metrics $\ip_{\mu,\nu}$ the tangent vector of
$\tilde \alpha$,
\begin{gather*}
  \tilde \alpha'(s) = \tilde \alpha(s) \cdot \bsmat \m 0 & -i\\ -i & \m 0\esmat,
\end{gather*}
is perpendicular to the vertical space at $\tilde \alpha(s)$,
which is spanned by the three vectors
\begin{align}
\begin{split}
\label{vertical}
  \xi_1(s) &:= \tfrac{d}{d\tau}\bigl( e^{i\tau} \cdot \tilde \alpha(s) \cdot
    \bsmat e^{\!-i\tau} & \m 0\\ 0 & \m 1\esmat \bigr)_{\vert \tau = 0}
    = \tilde \alpha(s) \cdot \bsmat -2i & \m 0\\ \m 0 & -i\esmat\\
  \xi_2(s) &:= \tfrac{d}{d\tau}\bigl( e^{j\tau} \cdot \tilde \alpha(s) \cdot
    \bsmat e^{\!-j\tau} & \m 0\\ 0 & \m 1\esmat \bigr)_{\vert \tau = 0}
    = \tilde \alpha(s) \cdot \bsmat j(\cos 2s - 1) & \m k \sin 2s\\ k \sin 2s & j\cos 2s\esmat\\
  \xi_3(s) &:= \tfrac{d}{d\tau}\bigl( e^{k\tau} \cdot \tilde \alpha(s) \cdot
    \bsmat e^{\!-k\tau} & \m 0\\ 0 & \m 1\esmat \bigr)_{\vert \tau = 0}
    = \tilde \alpha(s) \cdot \bsmat -k(\cos 2s + 1) & \m j \sin 2s\\
    \m j \sin 2s & -k \cos 2s\esmat. \qedhere
\end{split}
\end{align}
\end{proof}

\begin{lem}
The curve $\tilde \alpha$ is a geodesic for any of the metrics
$\ip_{\mu,\nu}$ on $\Syp(2)$.
\end{lem}

\begin{proof}
Since $\tilde \alpha$ is an integral curve of the left invariant vector field $v$ given by
$\bsmat \m 0 & -i\\ -i & \m 0\esmat \in \asyp(2)$ it suffices to compute
$\nabla_v v$ at the identity matrix. For an arbitrary left invariant vector field $w$ the
Kozul formula for the Levi-Civita connection yields
\begin{gather*}
  \langle \nabla_v v, w\rangle_{\mu,\nu}
  = - \langle v,[v,w]\rangle_{\mu,\nu}.
\end{gather*}
Using the special value of $v$ at the identity matrix and the
fact that $\ad_v$ is skew symmetric with respect to the biinvariant metric
$\ip_{\frac{1}{2},\frac{1}{2}}$ one gets
\begin{gather*}
  \langle \nabla_v v, w\rangle_{\mu,\nu}
  = -\langle v, [v, w]\rangle_{\mu,\nu}
  = -\langle v,  [v,w]\rangle_{\frac{1}{2},\frac{1}{2}}
  = -\langle [v,v], w\rangle_{\frac{1}{2},\frac{1}{2}} = 0
\end{gather*}
at the identity matrix.
\end{proof}

\begin{cor}
The curve $\alpha := \pi_{\Sigma^7}\circ \tilde \alpha$ is a geodesic in $\Sigma^7$,
and this geodesic is contained in the 5-sphere $\Sigma^5 \subset \Sigma^7$.
\end{cor}

\begin{lem}
The geodesic $\alpha$ in $\Sigma^5\subset \Sigma^7$ intersects all
$\bullet$-orbits perpendicularly.
\end{lem}

\begin{proof}
The tangent space to the $\bullet$-orbit through $\tilde \alpha(s)$ is spanned by
\begin{align}
\begin{split}
\label{orbitbasis}
  \hat v_0(s) &:=
  \tfrac{d}{d\theta} \bigl( \bsmat \cos\theta & -\sin\theta\\ \sin\theta & \m \cos\theta\esmat
      \cdot \tilde \alpha(s)\bigr)_{\vert\theta = 0}
    = \tilde \alpha(s)\cdot \bsmat i\sin 2s &\,  -\cos 2s \\ \cos 2s & \, -i\sin 2s\esmat,\\
  \hat v_1(s) &:=
 \tfrac{d}{d\tau} \bigl(\tilde \alpha(s) \cdot
   \bsmat 1 & 0\\ 0 & \m e^{-i\tau}\esmat \bigr)_{\vert\tau = 0}
    = \tilde \alpha(s) \cdot \bsmat 0 & \m 0 \\ 0 & -i\esmat,\\
  \hat v_2(s) &:=
 \tfrac{d}{d\tau} \bigl(\tilde \alpha(s) \cdot
   \bsmat 1 & 0\\ 0 & \m e^{-j\tau}\esmat \bigr)_{\vert\tau = 0}
    = \tilde \alpha(s)\cdot \bsmat 0 & \m 0 \\ 0 & -j\esmat,\\
  \hat v_3(s) &:=
 \tfrac{d}{d\tau} \bigl(\tilde \alpha(s) \cdot
   \bsmat 1 & 0\\ 0 & \m e^{-k\tau}\esmat \bigr)_{\vert\tau = 0}
    = \tilde \alpha(s)\cdot \bsmat 0 & \m 0 \\ 0 & -k\esmat.
\end{split}
\end{align}
All four vectors are perpendicular to the horizontal vector $\tilde \alpha'(s)$.
\end{proof}

The isotropy groups of the $\bullet$-action along the geodesic $\alpha$
are regular for $s\not\in\tfrac{\pi}{4}\cdot \Z$ and are denoted by $H$.
The singular isotropy groups at $s = 0$ and $s= \tfrac{\pi}{4}$ are denoted
by $K_-$ and $K_+$, respectively. Straightforward computations yield
\begin{lem}
\label{isotropy}
The isotropy groups along the normal geodesic $\alpha$ are given by
\begin{align*}
  H &= \bigl\{(\1,\pm 1), (-\1,\pm i), \bigl( \bsmat 1 & \m 0\\ 0 & -1\esmat, \pm j\bigr),
  \bigl( \bsmat -1 & \m 0\\ \m 0 & \m 1\esmat, \pm k\bigr)\bigr\}\\
  &\approx \Z_2 \times \Z_2,\\
   K_{-} &=  \{(\1, \pm e^{j\tau})\} \cup \{(-\1, \pm i e^{j\tau})\}
   \cup \bigl\{\bigl( \bsmat 1 & \m 0\\ 0 & -1\esmat, \pm e^{j\tau}\bigr)\bigr\}
   \cup \bigl\{\bigl( \bsmat -1 & \m 0\\ \m 0 & \m 1\esmat, \pm i e^{j\tau}\bigr)\bigr\}\\
   &\approx \Z_2 \times \OO(2),\\
  K_{+} &= \bigl\{ \bigl( \bsmat \cos \theta & -\sin \theta\\ \sin \theta & \m \cos \theta\esmat,
    \pm e^{-\frac{3}{2}i\theta} \bigr)\bigr\}
    \cup  \bigl\{ \bigl( \bsmat \cos \theta & -\sin \theta\\ \sin \theta & \m \cos \theta\esmat \cdot
    \bsmat 1 & \m 0\\ 0 & -1\esmat,
    \pm e^{-\frac{3}{2}i\theta}j \bigr)\bigr\}.\\
    & \approx \OO(2).
\end{align*}
\end{lem}
Note that $K_+$ is isomorphic to an $\OO(2)$ that projects surjectively onto
the $\OO(2)$ factor in the definition of the $\bullet$-action while $K_-$ is
isomorphic to $\Z_2\times \OO(2)$ where the $\OO(2)$-factor corresponds to
$\{(\1, \pm e^{j\tau})\} \cup \{(-\1, \pm i e^{j\tau})\}$, which is contained in the
identity component of the acting group $\OO(2)\times\SO(3)$.
The singular orbit at $s = 0$ is diffeomorphic to $\Sph^2\times_{\Z_2}\Sph^1$
and the singular orbit at $s = \tfrac{\pi}{4}$ is diffeomorphic to $\SO(3)$.

\medskip

We will now compute the induced metrics on the principal orbits.
This computation will be used later in this section to determine the isometry
group of the Gromoll-Meyer sphere and in section\,\ref{totgeod} where we
discuss the geometric properties of the metrics on $\Sigma^5$ and $\Sigma^6$.
We need to compute the inner products of four linearly independent Killing fields
along the normal geodesic $\alpha$ in $\Sigma^5$.
Such Killing fields $v_0(s),\ldots,v_3(s)$ are given by the horizontal parts
$\tilde v_0(s),\ldots,\tilde v_3(s)$ of the Killing fields
$\hat v_0(s),\ldots,\hat v_3(s)$ along $\tilde\alpha$ given in (\ref{orbitbasis}).
Straightforward computations using the orthogonal basis
$\xi_1(s)$, $\xi_2(s)$, $\xi_3(s)$
of the vertical space at $\tilde \alpha(s)$ given in (\ref{vertical}) show
\begin{align*}
  \tilde v_0(s) &= \tilde \alpha(s) \cdot
    \bigl( \tfrac{3\sin 2s}{4\mu+\nu} \bmat i \nu & 0\\
      0 & -2 i \mu \emat + \cos 2s \bmat 0 & -1\\ 1 & \m 0\emat \bigr),\\
  \tilde v_1(s) &= \tilde \alpha(s) \cdot
    \tfrac{2}{4\mu+\nu} \bmat i \nu & \, \m 0\\  0 & \,-2 i \mu \emat,\\
  \tilde v_2(s) &= \tilde \alpha(s) \cdot \Bigl( \bmat 0 & \m 0\\ 0 & -j\emat +
    \tfrac{\cos 2s}{\nu \cos^2 2s + 4 ( 1 - (1-\mu)\sin^2 s ) \sin^2 s}
    \bmat j(\cos 2s-1) & \;k \sin 2s\\ k \sin 2s & \;j\cos 2s \emat\Bigr),\\
  \tilde v_3(s) &=\tilde \alpha(s) \cdot \Bigl( \bmat 0 & \m 0\\ 0 & -k\emat +
    \tfrac{\cos 2s}{\nu \cos^2 2s + 4 ( 1 - (1-\mu)\cos^2 s ) \cos^2 s}
    \bmat k(1+\cos 2s) & \,-j \sin 2s\\  -j \sin 2s & \,\m k\cos 2s \emat\Bigr).
\end{align*}
The action of the principal isotropy group $H$ on these four Killing fields
along $\alpha$ is given by the matrices
\begin{gather*}
  \bmat 1 & \m 0 & \m 0 & \m 0\\ 0 & \m 1 & \m 0 & \m 0\\
    0 & \m 0 & \m 1 & \m 0\\ 0 & \m 0 & \m 0 & \m 1 \emat, \quad
  \bmat 1 & \m 0 & \m 0 & \m 0\\ 0 & \m 1 & \m 0 & \m 0\\
    0 & \m 0 & -1 & \m 0\\ 0 & \m 0 & \m 0 & -1 \emat, \quad
  \bmat -1 & \m 0 & \m 0 & \m 0\\ \m 0 & -1 & \m 0 & \m 0\\
    \m 0 & \m 0 & \m 1 & \m 0\\ \m 0 & \m 0 & \m 0 & -1 \emat, \quad
  \bmat -1 & \m 0 & \m 0 & \m 0\\ \m 0 & -1 & \m 0 & \m 0\\
    \m 0 & \m 0 & -1 & \m 0\\ \m 0 & \m 0 & \m 0 & \m 1 \emat.
\end{gather*}
The inner products $\langle \tilde v_j(s),\tilde v_k(s)\rangle$ of the Killing fields
are given by the matrix
\begin{gather}
\label{metricmatrix}
  \bmat
    a(s) & \,b(s)\, & 0 & 0\\
    b(s) & \frac{4\mu\nu}{4\mu+\nu} & 0 & 0\\
    0 & 0 & c(s) & 0\\
    0 & 0 & 0 & d(s)\emat
\end{gather}
where
\begin{align*}
  a(s) &= 1 - \bigl( 1-\tfrac{9\mu\nu}{4\mu+\nu} \bigr) \sin^2 2s,\\
  b(s) &= \tfrac{6\mu\nu}{4\mu+\nu}\,\sin 2s,\\
  c(s) &= \nu\,\tfrac{4 ( 1 - (1-\mu)\sin^2 s ) \sin^2 s}{\nu \cos^2 2s
    \,+\, 4 ( 1 - (1-\mu)\sin^2 s ) \sin^2 s}\,,\\
  d(s) &= \nu\,\tfrac{4 ( 1 - (1-\mu)\cos^2 s ) \cos^2 s}{\nu \cos^2 2s
    \,+\, 4 ( 1 - (1-\mu)\cos^2 s ) \cos^2 s}\,.
\end{align*}

This matrix description of the cohomogeneity one metrics $\ip_{\mu,\nu}$ on $\Sigma^5$
will be interpreted in section\,\ref{totgeod} in terms of totally geodesic submanifolds
$L^3$, $\Sigma^2$, and $\tilde \Sigma^2$, which intersect pairwise perpendicularly
in the normal geodesic $\alpha$. The upper left $2\times 2$-block of the matrix in
(\ref{metricmatrix}) describes cohomogeneity one metrics on the lens space
$L^3 \approx \Sph^3/\Z_3$ (the block becomes singular at $s \in \tfrac{\pi}{4} + \tfrac{\pi}{2}\Z$;
the smoothness at these times can best be seen by passing from $\tilde v_0$ to
$2 \tilde v_0 -3 \tilde v_1$ and from $s$ to $s + \frac{\pi}{4}$).
The numbers $c(s)$ and $d(s)$ describe cohomogeneity one metrics on the
$2$-spheres $\Sigma^2$ and $\tilde\Sigma^2$.
In section\,\ref{totgeod} we will see the following:
\begin{lem}
For all $\mu,\nu > 0$ the lens space $L^3 \approx \Sph^3/\Z_3$ is totally geodesic
in $\Sigma^5$ and $\Sigma^6_{\pm A}$ and inherits an intrinsically homogeneous metric.
For $\mu = 1$ and $\nu = \tfrac{1}{2}$ the lens space $L^3$ has constant curvature $1$.
\end{lem}
It is interesting to compare these metrics on $\Sigma^5$ (and hence on the exotic
projective space $\Sigma^5/\{\pm \1\}$\,) briefly to those that come from the Grove-Ziller
construction for cohomogeneity manifolds with codimension~2 singular orbits \cite{grove}.
Our metrics are analytic and there are always points with negative sectional curvature.
The Grove-Ziller metrics are merely smooth but have nonnegative sectional curvature;
on the lens space $L^3$ they induce a proper cohomogeneity one metric with planes
of zero sectional curvature over an open set of points.

\begin{lem}
\label{isomsigmafive}
The isometry group of $(\Sigma^5,\ip_{\mu,\nu})$ is the group
$\OO(2)\times \SO(3)$.
\end{lem}

\begin{proof}
It is obvious from the isotropy group computation above that the $\bullet$-action
of $G = \OO(2)\times \SO(3)$ on $\Sigma^5$ is effective.
The full isometry group $\tilde G$ of $\Sigma^5$ cannot act transitively
on $\Sigma^5$. Otherwise all fixed point sets of isometries would be homogeneous
which they are not (see section\,\ref{totgeod}). Hence, the geodesic $\alpha$
is perpendicular to all $\tilde G$-orbits. Let $\tilde H$ denote the common
principal isotropy group along the geodesic $\alpha$ and let $\tilde H(s)$
denote the group of isomorphisms of $\R^4$ that preserve the symmetric bilinear
form given by the matrix~(\ref{metricmatrix}). Clearly, $\tilde H$ is isomorphic to a
subgroup of the intersection of all $\tilde H(s)$. It is straightforward to see that the
intersection of all $\tilde H(s)$ is the group of order~$8$ generated
by the three involutions
\begin{gather*}
  \bmat -1 & \m 0 & \m 0 & \m 0\\ \m 0 & -1 & \m 0 & \m 0\\
    \m 0 & \m 0 & \m 1 & \m 0\\ \m 0 & \m 0 & \m 0 & \m 1 \emat, \quad
  \bmat 1 &  \m 0 & \m 0 & \m 0\\ 0 & \m 1 & \m 0 & \m 0\\
    0 & \m 0 & -1 & \m 0\\ 0 & \m 0 & \m 0 & \m 1 \emat, \quad
  \bmat 1 &  \m 0 & \m 0 & \m 0\\ 0 & \m 1 & \m 0 & \m 0\\
    0 & \m 0 & \m 1 & \m 0\\ 0 & \m 0 & \m 0 & -1 \emat.
\end{gather*}
If this entire group of order $8$ were the isotropy group along $\alpha$
then the curvature tensor $\langle R(\alpha',v_1)v_2,v_3\rangle_{\mu,\nu}$
of $\Sigma^5$ would vanish identically. However, computations show that
\begin{gather*}
  \langle R(\alpha',v_1)v_2,v_3 \rangle_{\mu,\nu} = r_{\mu,\nu}(\cos 2s) \sin 4s
\end{gather*}
where $r_{\mu,\nu}$ is a rational function with
$r_{\mu,\nu}(0) = -\tfrac{-8\mu\nu}{(1+\mu)(4\mu+\nu)}$.
This implies that $\langle R(\alpha',v_1)v_2,v_3\rangle_{\mu,\nu}$
does not vanish for $s$ close to but not equal to~$\tfrac{\pi}{4}$.
It follows that the isotropy group $H$ of the $\bullet$-action is the full principal
isotropy group $\tilde H$ of $(\Sigma^5,\ip_{\mu,\nu})$. Since the principal orbits
$G/H$ are connected, $G$ is the full subgroup of $\tilde G$ that
preserves the principal orbits. Moreover, $\tilde G$ is a finite extension of $G$ and
$G \subset \tilde G$ is a normal subgroup. If $G$ were a proper subgroup of $\tilde G$
then $\tilde G / G$ would act nontrivially on the orbit space and the Weyl group of the
cohomogeneity one action of $\tilde G$ on $\Sigma^5$ would be larger than that of
the action of $G$. This is impossible, as one can see from the isotropy groups
in Lemma\,\ref{isotropy}.
\end{proof}

\bigskip

\section{The identification of $\Sigma^5$ with the Brieskorn sphere $W^5_3$}
\label{briesid}
We will now construct a $\OO(2)\times \SO(3)$-equivariant diffeomorphism
between the sphere $\Sigma^5\subset \Sigma^7$ and the Brieskorn sphere $W^5_3$
given by the equations
\begin{gather*}
  \tfrac{8}{9}\,z_0^3 + z_1^2 + z_2^2 + z_3^2 \;=\; 0,\\
  \tfrac{4}{3}\,\abs{z_0}^2 + \abs{z_1}^2 + \abs{z_2}^2 + \abs{z_3}^2 \;=\; \tfrac{4}{9}
\end{gather*}
in $\C^4 = \C \oplus \C^3$. It is crucial for section\,\ref{briesparam} that we have
modified the coefficients compared to the standard definition of $W^5_3$.
The advantage of our choice is that there exists an explicit formula for a
unit speed geodesic in~$W^5_3$ that intersects all orbits of the action
\begin{align}
\begin{split}
\label{actionbries}
  \OO(2) \times \SO(3) \times W^5_3 &\to W^5_3, \\
  \bigl( \bmat \cos \theta & -\sin \theta\\ \sin\theta & \m \cos\theta \emat, A\bigr)
    \cdot (z_0, z) \;&=\; (e^{2i\theta}z_0, e^{3i\theta} A z),\\
  \bigl( \bmat 1 & \m 0\\ 0 & -1\emat , A\bigr)
    \cdot (z_0, z) \;&=\; ( \bar z_0, A \bar z)
\end{split}
\end{align}
with $z \in \C^3$ perpendicularly. This action on $W^5_3$ has first been considered
by Calabi (see Bredon's survey \cite{bredon}). In the literature, however, almost
exclusively the subaction of the identity component $\SO(2)\times \SO(3)$ is used.
The additional $\Z_2$ symmetry causes the fixed point set of the principal isotropy group
to be $1$-dimensional. Hence, we have preferred normal geodesics and hence
canonical identifications between $\Sigma^5$ and $W^5_3$ (see Lemma\,\ref{ngeod} below).

Consider the curve
\begin{gather*}
  \beta(s) = \Bigl(-\tfrac{1}{2}\cos 2s, \tfrac{1}{6} \bmat 0 \\ 3\cos s - \cos 3s\\
  3 i \sin s + i \sin 3s \emat \Bigr)
\end{gather*}
in $W^5_3 \subset \C \oplus \C^3$. It is straightforward to check that $\beta$
is parametrized by arc length and to compute the isotropy groups along $\beta$.

\begin{lem}
\label{isobries}
The isotropy groups $H$ at $\beta(s)$ for $s \not\in \tfrac{\pi}{4}\Z$ , $K_-$ at $\beta(0)$,
and $K_+$ at $\beta(\tfrac{\pi}{4})$ are given by
\begin{align*}
  H &= \Bigl\{(\1,\1), \Bigl(-\1,\bmat 1 & \m 0 & \m 0\\ 0 & -1 & \m 0\\ 0 & \m 0 & -1\emat\Bigr),
  \Bigl( \bmat 1 & \m 0\\ 0 & -1\emat,
  \bmat -1 & \m 0 & \m 0\\ \m 0 & \m 1 & \m 0\\ \m 0 & \m 0 & -1\emat \Bigr),
  \Bigl( \bmat -1 & \m 0\\ \m 0 & \m 1\emat,
  \bmat -1 & \m 0 & \m 0\\ \m 0 & -1 & \m 0\\ \m 0 & \m 0 & \m 1\emat \Bigr) \Bigr\}\\
  &\approx \Z_2 \times \Z_2,\\
  K_{-} &= 
    \Bigl\{\Bigl(\1, \bmat \ast & \m 0 & \m \ast\\ 0 & \m 1 & \m 0\\
      \ast & \m 0 & \m \ast \emat \Bigr)\Bigr\}
    \cup \Bigl\{\Bigl(-\1, \bmat \ast & \m 0 & \m \ast\\ 0 & -1 & \m 0\\
      \ast & \m 0 & \m\ast \emat \Bigr)\Bigr\}\\
   & \qquad \qquad\qquad\qquad\qquad\quad
   \cup \Bigl\{\Bigl(\bmat 1 & \m 0\\ 0 & -1\emat,
   \bmat \ast & \m 0 & \m \ast\\ 0 & \m 1 & \m 0\\ \ast & \m 0 & \m \ast \emat \Bigr)\Bigr\}
    \cup \Bigl\{\Bigl(\bmat -1 & \m 0\\ \m 0 & \m 1\emat,
   \bmat \ast & \m 0 & \m \ast\\ 0 & -1 & \m 0\\
      \ast & \m 0 & \m\ast \emat \Bigr)\Bigr\}\\
   &\approx \Z_2 \times \OO(2),\\
  K_{+} &= \Bigl\{ \Bigl( D(\theta),
    \bmat 1 & \m 0\\ 0 & \m D(-3\theta) \emat \Bigr)\Bigr\}
    \cup \Bigl\{ \Bigl( D(\theta) \cdot \bmat 1 & \m 0\\ 0 & -1\emat,
    \bmat 1 & \m 0\\ 0 & \m D(-3\theta) \emat
      \cdot \bmat -1 & \m 0 & \m 0\\ \m 0 & \m 1 & \m 0\\ \m 0 & \m 0 & -1\emat \Bigr)\Bigr\}.\\
    & \approx \OO(2),
\end{align*}
where $D(\theta) = \bmat \cos \theta & -\sin \theta\\ \sin \theta & \m \cos \theta\emat$
denotes the rotation in $\R^2$ with angle $\theta$.
\end{lem}

\begin{lem}
\label{ngeod}
The curve $\beta$ is a unit speed geodesic in $W^5_3$ that intersects all orbits of
the $\OO(2)\times\SO(3)$-action perpendicularly.
\end{lem}

\begin{proof}
The fixed point set of the principal isotropy group $H$ clearly contains a geodesic
that intersects all orbits perpendicularly. This fixed point set is given~by
\begin{gather*}
  \Imag z_0 = 0,\quad  z_1 = 0,\quad \Imag z_2 = 0,\quad \Real z_3 = 0.
\end{gather*}
It is easy to check that this fixed point set is one dimensional and that $\beta$
maps into the fixed point set of $H$.
\end{proof} 

In the following theorem it is supposed that $\SO(3) = \Sph^3/\{\pm 1\}$ is identified with
the matrix group $\SO(3)$ by the action of $\Sph^3$ on the imaginary quaternions by conjugation.
\begin{thm}
\label{equivdiff}
The map
\begin{gather*}
  \Sigma^5 \to W^5_3, \quad
  (A,\pm q)\bullet \alpha(s) \mapsto (A,\pm q)\cdot \beta(s)
\end{gather*}
is a well-defined $\OO(2)\times\SO(3)$-equivariant diffeomorphism.
\end{thm}

\begin{proof}
This follows from the isotropy groups in Lemma\,\ref{isotropy} and
Lemma\,\ref{isobries}
\end{proof}

\begin{cor}
\label{joinfin}
There is an $\OO(2)\times \SO(3)$-equivariant homeomorphism
\begin{gather*}
  \Sph^1 \ast W^5_3 \to \Sigma^7.
\end{gather*}
Here, $\OO(2)$ acts on $\Sph^1$ in the canonical~way.
\end{cor}
\begin{proof}
This follows directly from Corollary\,\ref{joingeod} and Theorem\,\ref{equivdiff}.
\end{proof}

\bigskip

\section{Free actions on the Gromoll-Meyer sphere}
In this section we classify all closed subgroups of $\OO(2)\times\SO(3)$ that
act freely on $\Sigma^7$ and determine the homotopy type of the orbit spaces.
Recall that $\OO(2)\times \SO(3)$ is the full isometry group of
$(\Sigma^7,\ip_{\mu,\nu})$ and $(\Sigma^5,\ip_{\mu,\nu})$
and that all elements in $\OO(2)\times\SO(3)$ that are not contained in
$\SO(2)\times\SO(3)$ have fixed points in $\Sigma^5 \subset \Sigma^7$
since they reverse the orientation of~$\Sigma^5$. In $\SO(2)\times\SO(3)$
it suffices to consider the elements of a maximal torus.

\begin{figure}
\begin{center}
\mbox{\scalebox{0.7}{\includegraphics{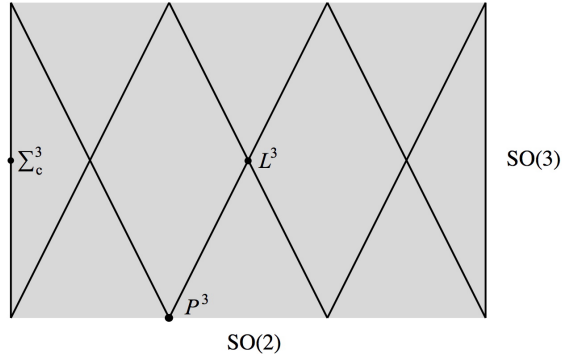}}}
\caption{The elements in the maximal torus of $\SO(2)\times\SO(3)$
with fixed points on $\Sigma^5$ are precisely illustrated by the black lines.}
\label{schrott}
\end{center}
\end{figure}

\begin{lem}
\label{sigmafix}
An element in $\SO(2)\times\SO(3)$ has a fixed point in $\Sigma^7$ if and only if
it has a fixed point in $\Sigma^5$.
\end{lem}

\begin{proof}
It follows from the isotropy group $K_-$ determined in Lemma\,\ref{isotropy}
that any isometry $(\1,\pm q)$ has fixed points in $\Sigma^5$. All other
elements of $\SO(2)\times\SO(3)$ are covered by Lemma\,\ref{nowiedfix}.
\end{proof}

\begin{lem}
An element in $\SO(2)\times\SO(3)$ has fixed points on $\Sigma^5$ if and only if
it is conjugate to an element of the subset of the maximal torus of
$\SO(2)\times\SO(3)$ illustrated in Figure\,\ref{schrott}.
\end{lem}

\begin{proof}
This follows with a few considerations from the computation of the isotropy
groups of the $\OO(2)\times\SO(3)$-action on $\Sigma^5$ in Lemma\,\ref{isotropy}.
\end{proof}

\begin{cor}
Every finite group that acts freely and isometrically on $(\Sigma^5,\ip_{\mu,\nu})$
and equivalently on $(\Sigma^7,\ip_{\mu,\nu})$ is cyclic.
\end{cor}

\begin{proof}
Let $G$ be a finite subgroup of $\SO(2) \times \SO(3)$. From Figure\,\ref{schrott}
we see that the kernel of the projection from $G$ to the $\SO(2)$
factor is trivial.
\end{proof}

\begin{cor}
All finite cyclic groups act freely and isometrically on $(\Sigma^5,\ip_{\mu,\nu})$
and hence also on the Gromoll-Meyer sphere $(\Sigma^7,\ip_{\mu,\nu})$.
\end{cor}

\begin{proof}
This is evident from extending the pattern of Figure\,\ref{schrott} periodically
to all of~$\R^2$.
\end{proof}

\begin{cor}
For every $m \in \N$ there are $7$-dimensional exotic homotopy lens spaces
with fundamental group $\Z_m$ and nonnegative sectional curvature.
For every even $m\in \N$ there are $5$-dimensional exotic homotopy lens
spaces with fundamental group $\Z_m$ and nonnegative sectional curvature.
\end{cor}

\begin{proof}
It is well-known that the quotient of a homotopy sphere by a cyclic group is
homotopy equivalent to a lens space (cf. \cite{browder}).
Since the Gromoll-Meyer sphere is not diffeomorphic to the standard sphere,
its quotients by finite cyclic groups cannot be diffeomorphic to lens spaces.
This completes the proof in the $7$-dimensional case. In the $5$-dimensional
case note that our metrics $\ip_{\mu,\nu}$ on $\Sigma^5$ do not have $K\ge 0$.
By the Grove-Ziller construction \cite{grove}, however, there are
$\OO(2)\times \SO(3)$-invariant metrics on $\Sigma^5 \approx W^5_3$
with $K \ge 0$. The quotient of $W^5_3$ by the Calabi involution is
homotopy equivalent but not diffeomorphic (and not homeomorphic) to $\RP^5$.
As is easily seen from Figure\,\ref{schrott} the $j\in \Z_{2j}$ acts by
the Calabi involution for all free $\Z_{2j}$-actions as above.
Hence, all corresponding quotients are not diffeomorphic to lens spaces.
\end{proof}

In the rest of this section we will determine to which lens spaces the orbit spaces
of the $\Z_m$-actions on $\Sigma^7$ and $\Sigma^5$ are homotopy equivalent.
Following an idea of Orlik \cite{orlik} we construct $\OO(2)\times \SO(3)$-equivariant
(continuous) branched coverings $\Sigma^5 \to \Sph^5$ and $\Sigma^7 \to \Sph^7$.
Using these branched coverings we obtain maps of degree $m l +1$
for some positive integer $l$ from $\Sigma^7/\Z_m$ and $\Sigma^5/\Z_m$
to standard lens spaces. By a theorem of Olum \cite{olum}, the existence
of such maps implies the existence of homotopy equivalences.

Let $D(\theta)$ denote the counterclockwise rotation in $\R^2$ by the angle $\theta$.
Then the subgroup $H_{m;p,q}$ of $\SO(2)\times \SO(3)$ generated by the element
\begin{gather*}
  \psi_{m;p,q} = \bigl(D(\tfrac{2\pi}{m}p), \pm e^{i\frac{\pi}{m}q} \bigr)
  = \Bigl(D(\tfrac{2\pi}{m}p), \bmat 1 & 0\\ 0 & D(\frac{2\pi}{m}q) \emat \Bigr)
\end{gather*}
acts freely on $\Sigma^7$ if and only if $p\neq 0$, $3p-q\neq 0$, $3p + q \neq 0$,
$m$ and $p$ are relatively prime, $m$ and $3p - q$ are relatively prime,
and $m$ and $3p+q$ are relatively prime. This is precisely what Figure\,\ref{schrott}
expresses graphically.

\begin{lem}
If $m$ is not divisible by $6$ then the quotient of $W^5_3$ by the
free action of $H_{m;p,q}$ is homotopy equivalent to the lens space
$L^5_m(p,3p+q,3p-q)$.
\end{lem}

\begin{proof}
Suppose first that $m$ is not divisible by $3$. The map
\begin{gather*}
  \varphi : W^5_3 \to \Sph^5, \quad
  (z_0,z_1,z_2,z_3) \mapsto \tfrac{1}{\sqrt{2(1-\abs{z_0}^2)}} (\sqrt{2}z_1,z_2+iz_3,z_3+iz_2)
\end{gather*}
is a $3:1$-covering branched along the singular orbit of the
$\SO(2)\times \SO(3)$-action on $W^5_3$ given by $z_0=0$.
If we define a $\Z_m$-action on $\Sph^5$ by
\begin{gather*}
  \Z_m \times \Sph^5 \to \Sph^5, \quad
  \bigl(j+m\Z,(z_1,z_2,z_3)) \mapsto
  \bigl( e^{i\frac{2\pi j}{m}3p} z_1, e^{i\frac{2\pi j}{m}(3p+q)} z_2,  e^{i\frac{2\pi j}{m}(3p-q)} z_3\bigr)
\end{gather*}
then $\varphi$ is $\Z_m$-equivariant. The orbit space of this $\Z_m$-action
on $\Sph^5$ is a lens space generally denoted by $L^5_m(3p,3p-q,3p+q)$.
Since $m$ is not divisible by $3$ there exists a positive integer $r$ such
that $3r \equiv 1 \mod m$.
Identify $\Sph^5$ homeomorphically with the join $\Sph^1\ast \Sph^3$ and consider the
continuous map $\rho: \Sph^5\to \Sph^5$ of degree $r$ induced by map
$\Sph^1\ast\Sph^3\to \Sph^1\ast \Sph^3$, $(\lambda, w) \mapsto (\lambda^r,w)$.
We now obtain the commutative diagram
\begin{gather}
\label{branchedcov}
\begin{CD}
W^5_3 @>{\varphi}>> \Sph^5 @>{\rho}>> \Sph^5\\
@VVV @VVV @VVV\\
W^5_3/H_{m,p,q} @>>> L^5_m(3p,3p-q,3p+q) @>>>
L^5_m(3rp,3p-q,3p+q)\\
@. @.  = L^5_m(p,3p-q,3p+q)
\end{CD}
\end{gather}
which includes a map of degree $3r$ between $W^5_3/H_{m,p,q}$
and $L^5_m(p,3p-q,3p+q)$. By Theorem\,4 of \cite{olum} these two
spaces are homotopy equivalent. In the case where $m$ is not
divisible by $2$ we can proceed similarly by exchanging the role
of $z_0$ and $z_1$ in the definition of $\varphi$.
\end{proof}

\begin{cor}
$W^5_3/H_{7;1,0}$ and $W^5_3/H_{7;1,1}$ are not homotopy
equivalent.
\end{cor}

\begin{proof}
This follows from the homotopy classification of lens spaces,
see \cite{olum}.
\end{proof}

\begin{cor}
If $m$ is not divisible by $6$ then the quotient of $\Sigma^7$ by the
free action of $H_{m;p,q}$ is homotopy equivalent to the lens space
$L^7_m(p,p,3p-q,3p+q)$.
\end{cor}

\begin{proof}
By Corollary\,\ref{joinfin}, $\Sigma^7$ is equivariantly homeomorphic to the
join $\Sph^1\ast W^5_3$. We obtain the statement by joining all spaces
in the diagram (\ref{branchedcov}) with $\Sph^1$ (note that $\psi_{m;p,q}$
acts on the circle $\Sigma^1\subset \Sigma^7$ as the rotation $D(\tfrac{2\pi}{m}p)$).
\end{proof}

The condition ``$m$ is not divisible by $6$'' seems to be a technical artifact.

\begin{cor}
$\Sigma^7/H_{5;1,0}$ and $\Sigma^7/H_{5;1,1}$ are not homotopy
equivalent.
\end{cor}

\bigskip

\section{Fixed point sets of isometries}
\label{totgeod}
Recall from section\,\ref{isomgroup} that there are three types of
isometries of $\Sigma^7$: Isometries that do not have fixed points
in $\Sigma^1$ (type I), isometries that fix precisely two points in $\Sigma^1$
(type II), and isometries that fix all points in $\Sigma^1$ (type III).

Isometries of type III are of the form $(\1,\pm q)$ with $q \in \Sph^3$, $q\neq \pm 1$.
They correspond up to conjugation to the black vertical line in Figure\,\ref{schrott}.
From Lemma\,\ref{isotropy} we see that the fixed point set in $\Sigma^5$
is a circle that is located in the singular orbit through $\alpha(0)$.
Thus by Lemma\,\ref{fullwiedfix} the fixed point set in $\Sigma^7$ is the join of
$\Sigma^1$ and this circle and hence diffeomorphic to $\Sph^3$.
Although not all elements of the form $(\1,\pm q)$ are conjugate to $(\1,\pm i)$,
all their fixed point sets are congruent to the fixed point set of $(\1,\pm i)$,
which will be denoted by $\smc$.

\begin{lem}
The fixed point set $\smc$ of $(\1,\pm i)$ on $\Sigma^7$ is isometric to
a $3$-sphere equipped with a Berger metric where the horizontal geodesics
have length $2\pi$ and the Hopf circles have length $2\pi\sqrt{\mu}$.
\end{lem}
\begin{proof}
With the structural information above it is immediate that $\smc = \pi_{\Sigma^7}(\U(2))$,
which is isometric to the homogeneous space $\U(2)/\U(1)$ where $\U(1)$ is
embedded into the right lower corner and $\U(2)$ is equipped with the
metric $\ip_{\mu,\nu}$.
\end{proof}

\begin{cor}
\label{eqhombries}
There is no $\SO(3)$-equivariant homeomorphism between $\Sigma^7$
and any of the Brieskorn spheres $W^7_{6j-1,3}$.
\end{cor}
\begin{proof}
In $W^7_{6j-1,3}$ the fixed point set of $\pm i = \diag(1,-1,-1)\in \SO(3)$
is the integral homology sphere $W^3_{6j-1,3} = W^3_{6j-1,3,2}$.
For $j = 1$ this space is diffeomorphic to Poincare dodecahedral space
and for $j > 1$ the universal cover of $W^3_{6j-1,3,2}$ is
$\widetilde\SL(2,\R)$ (see \cite{milnor}).
\end{proof}
Note that this last argument also gives a simple reason for why there are
no $\SO(3)$-invariant Riemannian metrics on $W^7_{6j-1,3}$ with $K > 0$
for $j>1$.

\begin{lem}
The circle $\Sigma^1$ is a closed geodesic for all $\SO(3)$-invariant
Riemannian metrics on $\Sigma^7$, in particular for all metrics $\ip_{\mu,\nu}$.
\end{lem}
\begin{proof}
It is easy to check that $\Sigma^1$ is the intersection of all the fixed point sets
$(\Sigma^7)^{(1,\pm q)}$, i.e., the common fixed point set of~$\SO(3)$.
\end{proof}

Isometries of type I are contained in $\SO(2)\times \SO(3)$.
Those with fixed points correspond up to conjugation to the skew lines in Figure\,\ref{schrott}.
It suffices to consider the left half of Figure\,\ref{schrott} since one can conjugate any
isometry by $\bigl( \bmat 1 & \m 0\\ 0 & -1\emat, \1\bigr)\in \OO(2)\times\SO(3)$.
The isometry $(-\1,\pm i)$ corresponds to the midpoint of the torus in Figure\,\ref{schrott}.
Its fixed point set $L^3$ is diffeomorphic to a lens space $\Sph^3/\Z_3$.
The fixed point sets of the remaining isometries of type I are all contained
in the singular orbit through $\alpha(\tfrac{\pi}{4})$ by Lemma\,\ref{isotropy}.
Up to conjugation only the fixed point set $P^3$ of the isometry
$\Bigl(\Bsmat \cos\frac{2\pi}{3} & - \sin \frac{2\pi}{3}\\
  \sin \frac{2\pi}{3} & \m \cos \frac{2\pi}{3}\Esmat, \1\Bigr)$
is more than $1$-dimensional.

\begin{lem}
The fixed point set $P^3$ of the isometry $\Bigl(\Bsmat \cos\frac{2\pi}{3} &
- \sin \frac{2\pi}{3}\\ \sin \frac{2\pi}{3} & \m \cos \frac{2\pi}{3}\Esmat, \1\Bigr)$
on $\Sigma^7$ is isometric to $\RP^3$ covered by a Berger $\Sph^3$
whose horizontal geodesics have the length $2\pi\sqrt{\nu}$
and whose Hopf circles have length $2\pi \sqrt{\tfrac{4\mu\nu}{4\mu+\nu}}$.
\end{lem}
\begin{proof}
It is immediate from the isotropy groups along the normal geodesic $\alpha$
in Lemma\,\ref{isotropy} that $P^3$ is precisely the $\OO(2)\times\SO(3)$-orbit
through $\alpha(\tfrac{\pi}{4})$. The subgroup $\SO(3)$ acts simply transitively
on this orbit and the induced metric can be obtained from (\ref{metricmatrix}).
\end{proof}

\begin{lem}
The fixed point set $L^3$ of the isometry $(-\1,\pm i)$ on $\Sigma^7$
is diffeomorphic to the lens space $\Sph^3/\Z_3$ and totally geodesic in
$\Sigma^5$, $\Sigma^6_{\pm A}$, and $\Sigma^7$. Moreover, $L^3$ is isometrically
covered by a Berger metric on $\Sph^3$ where the horizontal geodesics have length
$2\pi$ and the Hopf circles have length $2\pi\sqrt{\tfrac{9\mu\nu}{4\mu+\nu}}$.
In particular, the extremal values of the sectional curvature of $L^3$ at any point
are $\tfrac{9\mu\nu}{4\mu+\nu}$ and $4 - \tfrac{27\mu\nu}{4\mu+\nu}$.
\end{lem}

\begin{proof}
In section \ref{briesid} it is shown that $\Sigma^5$ and the
Brieskorn sphere $W^5_3$ are equivariantly diffeomorphic.
In $W^5_3$ the corresponding fixed point set is $W^3_3$ which is
well-known to be diffeomorphic to $\Sph^3/\Z_3$ (see e.g.\ \cite{mayer}).
In $\Sigma^5$ there exists a direct argument that allows a simple
curvature computation: Straightforward computations show that
the horizontal lift of $T_{\alpha(s)} L$ at $\tilde\alpha(s)$ is spanned by
the horizontal vectors $\tilde \alpha'(s)$, $\tilde v_0(s)$, $\tilde v_1(s)$.
Thus $L$ is 3-dimensional. Let $\U(2)$ be the centralizer of $\bmat i & 0\\ 0 & i\emat$
in $\Syp(2)$. It is straightforward to see that $\pi_{\Sigma^7}(j\U(2))$ embeds
into the fixed point set. Now $\pi_{\Sigma^7}(j\U(2))$ is isometric to the quotient
of $j\U(2)\subset \Syp(2)$ by the $\U(1)$-action
$(\lambda,jA) \mapsto \lambda jA \bmat \bar \lambda & 0\\ 0 & 1\emat$
where $\U(2)$ carries the metric $\ip_{\mu,\nu}$ induced from $\Syp(2)$.
Since $i$ and $j$ anticommute this quotient is isometric the homogeneous space
$\U(2)/\bigl\{\bsmat \bar \lambda^2 & 0\\ 0 & \bar \lambda\esmat \bigr\}$
and hence diffeomorphic to $\Sph^3/\Z_3$. The vector
$\tilde\alpha'(0) = j\bmat \m 0 & -i\\ -i & \m 0\emat$
is horizontal with respect to the fibration
\begin{gather*}
  \U(2)/\bigl\{\bsmat \bar \lambda^2 & 0\\ 0 & \bar \lambda\esmat \bigr\}
  \to \U(2)/(\U(1)\times\U(1)) = \CP^1
\end{gather*}
and the curve $\exp(t\tilde\alpha'(0))$ closes first after length $2\pi$ in $\U(2)$.
The vector
$\tilde v_1(0) = j \tfrac{2}{4\mu+\nu} \bmat i\nu & 0\\ 0 & -2i\mu\emat$
is vertical. It is easy to check that the curve $\exp(t\tilde v_1(0))$ in $\U(2)$ meets
the circle $\bigl\{\bsmat \bar \lambda^2 & 0\\ 0 & \bar \lambda\esmat \bigr\}$
first at time $T = \pi$. Hence, the length of the vertical circle in the lens space
is $\pi\abs{\tilde v_1(0)} = 2\pi \sqrt{\tfrac{\mu\nu}{4\mu+\nu}}$. In the universal
cover $\Sph^3$ the length of this Hopf circle is three times as long.
\end{proof}

We finally deal with isometries of type II, i.e., isometries that fix precisely two points
in $\Sigma^1$. It is clear that these isometries are not contained in $\SO(2)\times\OO(3)$.
By Lemma\,\ref{isotropy} and Lemma\,\ref{twofixwied} they are conjugate to
$\bigl( \bmat 1 & \m 0\\ 0 & -1\emat, \pm 1\bigr)$ or
$\bigl( \bmat 1 & \m 0\\ 0 & -1\emat, \pm j\bigr)$
if the dimension of the fixed point set is $> 1$. From Lemma\,\ref{twofixwied}
it is also clear that the fixed point sets of these isometries are suspensions of
subspheres of $\Sigma^5$ from the two points $\Sph^3\star (\pm \1)$.

\begin{lem}
The fixed point set $\sme$ of the isometry $\bigl( \bmat 1 & \m 0\\ 0 & -1\emat, \pm 1\bigr)$
is diffeomorphic to $\Sph^3$. The induced metric on
$\sme \smallsetminus (\Sph^3\star (\pm \1))$ is isometric to the metric
\begin{gather*}
  \mu \bigl( dt^2 + \tfrac{\nu \sin^2 t}{\nu + 4 \mu \sin^2 t} \, g^{\Sph^2}_{\mathrm{can}}\bigr)
\end{gather*}
on $[0,\pi] \times \Sph^2$. Hence, the sectional curvatures vary
between $\tfrac{\nu}{\mu(4\mu+\nu)}$ and $\tfrac{12\mu+\nu}{\mu\nu}$.
\end{lem}

\begin{proof}
By Lemma\,\ref{isotropy} the fixed point set is the suspension of homogeneous
$2$-spheres from the two points $\Sph^3\star(\pm \1)$. This suspension is
given by the geodesics $\pi_{\Sigma^7}\circ\tilde\gamma_{\bmat p \\ 0 \emat}$ in (\ref{lift}).
It is straightforward to compute the diameter of the $\SO(3)$-orbits through
$\pi_{\Sigma^7}\circ\gamma_{\bmat p \\ 0 \emat}(t)$.
\end{proof}

\begin{lem}
\label{sigmatwo}
The fixed point set $\Sigma^2$ of the isometry
$\bigl( \bmat 1 & \m 0\\ 0 & -1\emat, \pm j\bigr)$ on $\Sigma^5$ and on $\Sigma^6_{\pm \1}$
is isometric to a $2$-sphere equipped with the metric $ds^2 + \tfrac{1}{4}c(s) d\phi^2$.
Here, $c(s)$ is the function $[0,\pi] \to \R$ defined in {\rm (\ref{metricmatrix})}.
The sectional curvature $K$ of $\Sigma^2$ satisfies
\begin{gather*}
  K_{\vert s = 0} = \tfrac{12}{\nu} - 8 - 3\mu, \quad
  K_{\vert s = \frac{\pi}{4}} = \tfrac{4\nu}{1+\mu}, \quad
  K_{\vert s = \frac{\pi}{2}} = -\tfrac{\nu (1+2\mu)}{\mu (4\mu+\nu)}.
\end{gather*}
\end{lem}

\begin{proof}
The manifold structure of the fixed point set in $\Sigma^5$ can best be determined by passing
from $\Sigma^5$ to the Euclidean sphere $\Sph^5$ with the nonlinear action obtained
in section\,\ref{nlEucl}. On $\Sph^5\subset \Imag\H\times\H$ it is straightforward
to check that the transformation $\bigl( \bmat 1 & \m 0\\ 0 & -1\emat, \pm j\bigr)$
fixes precisely the 2-sphere that consists of all unit vectors of the form
$\bsmat p\\w\esmat$ with $p\in j\R$ and $w \in \spann_{\R}\{i,k\}$.
The metric on the fixed point set $\Sigma^2\subset \Sigma^5$, however, has to be
determined in $\Sigma^5$. It is easy to see that $\Sigma^2$ contains the normal
geodesic $\alpha$ and that the tangent space to $\Sigma^2$ at $\alpha(s)$
is spanned by $\alpha'(s)$ and the Killing field $v_2(s)$ if $s \not \in \pi\Z$.
A straightforward computation shows that the circle which corresponds to $v_2$
and acts effectively on $\Sigma^2$ inherits the length $\pi \cdot \sqrt{c(s)}$
at time $s$. This completes the computation of the induced metric.
The curvature computations are straightforward. Finally, Lemma\,\ref{nowiedfix}
assures that the fixed point set of the isometry on $\Sigma^6_{\pm \1}$
is contained in $\Sigma^5$.
\end{proof}

\begin{cor}
\label{sigmathree}
The fixed point set $\sma$ of the isometry
$\bigl( \bmat 1 & \m 0\\ 0 & -1\emat, \pm j\bigr)$ on $\Sigma^7$
is diffeomorphic to $\Sph^3$.
\end{cor}

In order to describe the metric that $\sma$ inherits from $\Sigma^7$ it is useful
to note that the horizontal lift of the tangent space $T_{\Sph^3\star\1}\sma$ at $\1$
is spanned by the three vectors
\begin{gather*}
  \bmat j & 0\\ 0 & 0\emat, \quad
  \bmat 0 & i\\ i & 0\emat, \quad
  \bmat 0 & k\\ k & 0 \emat.
\end{gather*}
Hence, $\sma$ can be parametrized by the horizontal geodesics
$\tilde \gamma_{\bmat p\\ w\emat}$ given in (\ref{lift}) with
\begin{gather*}
  p = j \cos \theta \quad \text{and}\quad w = i \sin\theta \cos \phi + k \sin\theta \sin\phi
\end{gather*}
where $t \in [0,\pi]$, $\theta \in [0,\pi]$, and $\phi\in [0,2\pi]$.
Thus, $\sma$ corresponds to a maximal choice of anticommuting
$p$ and $w$. In the coordinates $(t,\theta,\phi)$ the metric on $\sma$ is given by
\begin{align*}
  g_{11} \,=\,& 1 - \tfrac{1}{D} (1-\mu) \bigl( 4 \sin^2 t \,\sin^2 \theta +
    \nu(1 - 2 \sin^2 t \, \sin^2 \theta)^2 \bigr) \cos^2 \theta, \\
  g_{22} \,=\,& \sin^2 t + \tfrac{1}{D} \sin^2 t \, \sin^2 \theta
    \Bigl( \nu \sin^2 \theta \, (2t-\sin 2t)^2 - (1-\mu) \,\cdot \\
   & \qquad\qquad \cdot \bigl( (\nu + 4\sin^2 t\, \sin^2 \theta) \cos^2 t
   + 2t\, \nu\,\sin^2\theta \, (2t \sin^2 t \, \sin^2 \theta - \sin 2t)\bigr)\Bigr),\\
 g_{33} \,=\,& \tfrac{1}{D} \, \nu \sin^2 t \,\sin^2 \theta \bigl( 1 - (1-\mu)\sin^2 t \,\sin^2 \theta\bigr),\\
 g_{23} \,=\,& \tfrac{1}{D} \, \nu \sin^2 t \,\sin^3 \theta \bigl( 2t - \sin 2t +
   \tfrac{1-\mu}{2} (\sin 2t - 4t \sin^2 t \, \sin^2 \theta ) \bigr),\\
 g_{13} \,=\,& -\tfrac{1}{D} \, \nu (1-\mu) \sin^2 t\, \sin^2 \theta\,
   \cos\theta \, (1-2\sin^2 t \, \sin^2\theta),\\
 g_{12} \,=\,& \tfrac{1-\mu}{4 D} \sin 2\theta \bigl( 4 \sin 2t \, \sin^2 t\,\sin^2 \theta \\
   & \qquad\qquad\quad - \nu (1-2\sin^2 t\,\sin^2\theta) (4t \sin^2t\,\sin^2\theta - \sin 2t)\bigr)
\end{align*}
where
\begin{gather*}
  D = 4 \bigl( 1 - (1-\mu) \sin^2 t \, \sin^2 \theta \bigr) \sin^2 t \,\sin^2 \theta
    + \nu(1 - 2 \sin^2 t \, \sin^2 \theta)^2.
\end{gather*}
This specializes for $\mu = 1$ to
\begin{gather*}
 g =  \Bsmat 1 & 0 & 0\\ 0 & \,\sin^2 t\, & 0\\ 0 & 0 & 0\Esmat
  + \tfrac{\nu \sin^2 t \sin^2 \theta}{4\sin^2 t \sin^2 \theta
    + \nu(1 - 2 \sin^2 t \sin^2 \theta)^2}
      \,\bbsmat 0\\ (2t-\sin 2t)\sin\theta \\ 1\Esmat \cdot
      \bbsmat 0\\ (2t-\sin 2t)\sin\theta \\ 1\Esmat^{\text{tr}}.
\end{gather*}

Note that $\sma$ is invariant under the isometry of $\Sigma^7$ induced by
$-\1 \in \Syp(2)$. In our coordinates, this isometry is given by
$(t,\theta,\phi) \mapsto (\pi - t, \theta + \pi, \phi - 2\pi \cos \theta)$.
This coordinate change allows us to glue $\sma$ from two disks
equipped with $g$.

Although the metric $g$ is of cohomogeneity~$2$ the curvature
of $g$ behaves like the action on $\sma$ were of cohomogeneity~$1$:
The orbit space of the natural $\SO(2)$-action on $\sma$
(in our coordinates given by translation in $\phi$) can easily be
shown to be the hemisphere of constant curvature $1$ for $\mu = 1$.
In our coordinates this hemisphere is given by $0 \le t \le \pi$, $0 \le \theta \le \pi$
and represented by geodesics from a point in the boundary.
It is reasonable to switch to polar coordinates, i.e., to make the coordinate change
\begin{gather*}
  \bmat \cos t \\ \cos \theta \sin t \\ \sin \theta \sin t \emat
  = \bmat \sin \omega \cos \psi \\ \sin \omega \sin \psi\\ \cos \omega \emat
\end{gather*}
with $0 \le \omega \le \tfrac{\pi}{2}$ and $0 \le \psi \le 2\pi$.
\begin{lem}
The orbit space of the natural $\SO(2)$-action on $(\sma,\ip_{\mu,\nu})$
is a hemisphere which inherits a rotationally invariant metric with curvature
\begin{gather*}
  \mu \, \tfrac{1+ 2(1-\mu)\cos^2 \omega}{(1-(1-\mu)\cos^2 \omega)^2} \, .
\end{gather*}
In particular, the curvature is constant if and only if $\mu=1$.
\end{lem}
The eigenvalues of the Einstein tensor (i.e., the critical values of the
sectional curvature) of $\sma$ also turn out to be independent of $\psi$.
The metric $g$ itself, however, does not improve in the coordinates $(\omega,\psi,\phi)$
nor does the curvature computation become simpler.
The next lemma gives some more detailed curvature information.
\begin{lem}
The scalar curvature of $(\sma,\ip_{1,\nu})$ is given by
\begin{gather*}
  \tfrac{4(-12+4\nu+9\nu^2 + 2(21\nu-8)\cos 2\omega + (9\nu^2 + 16\nu -4)\cos 4\omega
  + 2\nu \cos 6\omega)}{(4+\nu+4\cos 2\omega + \nu \cos 4\omega)^2}
\end{gather*}
For $\mu,\nu \le 1$ the minimum of the sectional curvature is given by
\begin{gather*}
  \min K = \min \bigl\{ \tfrac{\mu\nu}{4\mu+\nu}, \tfrac{12-8(\mu +\nu) - 3\mu\nu}{4\mu+\nu} \bigr\}.
\end{gather*}
In the Gromoll-Meyer case $\mu = \nu = \frac{1}{2}$, $\sma$
inherits a metric with $\tfrac{\min K}{\max K} = \tfrac{1}{145}$.
\end{lem}

Note that by construction and Lemma\,\ref{twofixwied}, $\sma$ is totally geodesic
in $\Sigma^7$ and in $\Sigma^6_{\pm A_0}$ with $A_0 = \bmat 0 & -1\\ 1 & \m 0\emat$,
and $\Sigma^2$ is totally geodesic in $\Sigma^5$ and in $\Sigma^6_{\pm \1}$.

\begin{cor}
There is a point in $\Sigma^2$ which has negative curvature for
all the metrics $\ip_{\mu,\nu}$ on $\Sigma^7$. Moreover, $\Sigma^5$ and
$\Sigma^6_{\pm \1}$ are not totally geodesic in $\Sigma^7$ for any of these metrics.
\end{cor}

\begin{proof}
The intrinsic sectional curvature of $\Sigma^2$ at $\alpha(\tfrac{\pi}{2})$
is $-\tfrac{\nu(1+2\mu)}{\mu(4\mu + \nu)} < 0$ by Lemma\,\ref{sigmatwo}.
The point $\alpha(\tfrac{\pi}{2})\in \Sigma^2\subset \sma$ corresponds to the
coordinates $t=\theta=\phi=\tfrac{\pi}{2}$ on $\sma$. The extrinsic sectional
curvature of the tangent space of $\Sigma^2$ at this point can be computed to
$\tfrac{\mu\nu}{4\mu+\nu} > 0$.
\end{proof}

We would like to add some comments on these fixed point sets:
First, the spheres $\smc$, $\sme$, and $\sma$ are indexed according to
their intrinsic cohomogeneity.
Second, the fixed point sets $\smc$, $L^3$, and $\sma$ yield necessary conditions
for $(\Sigma^7,\ip_{\mu,\nu})$ to have nonnegative sectional curvature:
That $\smc$ inherits $K\ge 0$ implies $\mu \le \tfrac{4}{3}$,
that $L^3$ inherits $K\ge 0$ implies $4(4\mu + \nu) - 27 \mu\nu \ge 0$,
and that $\sma$ inherits $K\ge 0$ implies $12 - 8 (\mu+\nu) - 3\mu\nu \ge 0$.
The last inequality is the most restricting one.

Of particular interest is the question of whether the nice behaviour of geodesics
on $\Sigma^7$ for $\mu = 1$ can be combined with nonnegative sectional curvature.
In this case the inequalities above show that necessarily $\nu \le \tfrac{4}{11}$
(this is precisely the inequality that guarantees that $\sma$ inherits $K\ge 0$).
For $\mu = 1$ and any $\nu > 0$ there are always some negative sectional
curvatures on $(\Syp(2),\ip_{\mu,\nu})$. The question whether these disappear
for small $\nu > 0$ when going down to $\Sigma^7$ seems to be subtle.

A distinguished metric on $\Sigma^7$ is $\ip_{1,\frac{1}{2}}$. In this case
$\smc$ and $L^3$ both have constant curvature $1$.

Finally, we would like to point out that not all totally geodesic submanifolds
of $\Sigma^7$ are fixed point sets of isometries:

\begin{lem}
For any metric $\ip_{\mu,\nu}$ the rectangular $2$-torus $T^2$ in $\Syp(2)$
para\-metrized by
\begin{gather*}
  \tfrac{1}{\sqrt{2}} \Bsmat 1\phantom{^2} & i\\ i\phantom{^2} & 1\Esmat \cdot
    \bmat e^{i\alpha} & 0\\ 0 & e^{j\beta}\emat \,=\, \tfrac{1}{\sqrt{2}}
    \bmat e^{i\alpha} & i e^{j\beta}\\ ie^{i\alpha} & e^{j\beta} \emat
\end{gather*}
with $\alpha,\beta \in \R$ is totally geodesic and horizontal with respect to the
submersion $\pi_{\Sigma^7}: \Syp(2) \to \Sigma^7$. Its image is a totally geodesic
rectangular $2$-torus in $\Sigma^7$ covered twice by $T^2$.
\end{lem}

\begin{proof}
Consider the subgroup $G$ of $\Syp(1)\times \Syp(1) \subset \Syp(2)$
generated by the two elements $\bmat i & 0\\ 0 & 1\emat$ and $\bmat 1 & 0\\ 0 & j\emat$.
This group $G$ acts by conjugation isometrically on $(\Syp(2),\ip_{\mu,\nu})$.
The rectangular torus
\begin{gather*}
  \Bigl\{ \bmat e^{i\alpha} & 0\\ 0 & e^{j\beta}\emat \,\Big\vert\; \alpha,\beta\in \R\Bigr\}
  \subset \Syp(2)
\end{gather*}
is the common fixed point set of $G$ and hence totally geodesic.
Hence, its left translated copy $T^2$ is totally geodesic, too. It is straightforward to show that
$T^2$ is horizontal and that $\pi_{\Sigma^7}$ restricted to $T^2$ induces an
embedding of $T^2/\bsmat 1 & 0\\ 0 & \pm 1\esmat$ into~$\Sigma^7$.
\end{proof}

This torus was already implicitly contained in \cite{meyer} and is also
listed in \cite{wilhelm}. The fundamental difference between the standard action $\bullet$
of $\Sph^3$ on $\Syp(2)$ and the Gromoll-Meyer action $\star$ appears here
very clearly: The torus $T^2$ is horizontal for the $\star$-action while only an
$\Sph^1$-factor is horizontal for the $\bullet$-action.

\bigskip

\section{An explicit parametrization of two Brieskorn spheres}
\label{briesparam}
In this section we present an explicit formula for two diffeomorphisms
between Euclidean spheres and Brieskorn spheres. The coefficients
in this formula are rational functions of the coordinates of the sphere.
They are simple enough that the entire formula fits into a few lines but
complicated enough that they could never be guessed.
The formula was obtained by combining the geodesic parametrization of
$\Sigma^5\subset \Sigma^7$ and the cohomogeneity one diffeomorphism
between $\Sigma^5$ and the Brieskorn sphere $W^5_3$.
The steps of the computations behind this approach will be explained at the
end of this section. The properties of the final formula, however, can be verified
straight, which shows that the formula is also valid in dimension~$13$
where no geometric derivation is possible so far.

Analogously to the previous section the Brieskorn sphere $W^{2n-1}_3$
is defined by the equations
\begin{gather*}
  \tfrac{8}{9}\,z_0^3 \,+\, z_1^2 \,+\, z_2^2 \,+ \ldots +\, z_n^2 \;=\; 0,\\
  \tfrac{4}{3}\,\abs{z_0}^2 + \abs{z_1}^2 + \abs{z_2}^2 + \ldots + \abs{z_n}^2
  \;=\; \tfrac{4}{9}
\end{gather*}
for $(z_0,z)\in \C \oplus \C^n$. For odd $n \ge 3$ the Brieskorn sphere $W^{2n-1}_3$
is diffeomorphic to the Kervaire sphere (see e.g. \cite{mayer}). By a result of
Brouwder, $W^{2n-1}_3$ can hence only be diffeomorphic to $\Sph^{2n-1}$
if $n = 2^{m}-1$. Up to now, this is known to hold for $n \in \{3,7,15,31\}$.
For $n = 3,7$ the classification theorems of J\"anich and the Hsiang brothers show
that there exist $\SO(3)$-equivariant diffeomorphisms $\Sph^5 \to W^5_3$
and $\Gtwo$-equivariant diffeomorphisms $\Sph^{13}\to W^{13}_3$.
We will construct the first explicit formulas for such diffeomorphisms here.

We decompose $z_0$ and $z$ into their real and imaginary parts, i.e.,
we set $z_0 = x_0 + i y_0$ and $z = x + iy$.
This leads to the equivalent definition of the Brieskorn sphere $W^{2n-1}_3$
by the three real equations
\begin{align}
\begin{split}
\label{realbrieskorn}
  \abs{x}^2 \,&=\, \tfrac{2}{9} \,(1 - 2 x_0^3 + 6 x_0 y_0^2 - 3 x_0^2 - 3 y_0^2),\\
  \abs{y}^2 \,&=\, \tfrac{2}{9} \,(1 + 2 x_0^3 - 6 x_0 y_0^2 - 3 x_0^2 - 3 y_0^2),\\
  \langle x,y \rangle \,&=\, \tfrac{4}{9}\,y_0\, (y_0^2 - 3x_0^2)
\end{split}
\end{align}
for $x_0,y_0 \in \R$ and $x,y\in \R^n$. The natural $\SO(n)$-action on $W^{2n-1}_3$
multiplies $x$ and $y$ by a matrix $A\in \SO(n)$ and leaves $x_0$ and $y_0$
unchanged. Analogously to \cite[pages 31--32]{mayer} one can show that
the orbit space of this action can be identified with the disc
$D^2 = \{ \lambda \in \C\,\big\vert\, \abs{\lambda} \le 1\}$
by the projection map $W^{2n-1}_3 \to D^2$, $(z_0,z) \mapsto 2 z_0$.
For $n = 7$ the action of $\Gtwo \subset \SO(7)$ has the same orbits
as the $\SO(7)$-action.

\smallskip

We now parametrize the standard spheres $\Sph^5$ and $\Sph^{13}$
by two vectors $p,w \in \R^3$ (resp.\ $\R^7$) with $\abs{p}^2+\abs{w}^2 = 1$
and set
\begin{align}
\begin{split}
\label{identbrieskorn}
  x_0 \,=\,& \tfrac{1}{2} (\abs{w}^2 - \abs{p}^2),\\
  y_0 \,=\,& - \langle p,w\rangle,\\
  x \,=\,& \tfrac{1}{3(1+\abs{p}^2)^2}
    \Bigl( \bigl( (3 - 2 \abs{p}^2) \, (1+\abs{p}^2)^2
    \,-\, 4 (1 - \abs{p}^2) \langle w,p \rangle^2 \bigl) \, p\\[-1ex]
    & \qquad\qquad\quad
      - 2 \bigl( 3 + 8 \abs{p}^2 + \abs{p}^4 - 4 \langle w,p \rangle^2 \bigr)
      \langle p,w \rangle \, w
    \,-\, 8 \abs{p}^2 \langle p,w \rangle \, p\times w \Bigr),\\
  y \,=\,& \tfrac{1}{3(1+\abs{p}^2)^2}
    \Bigl( \bigl( - ( 1+ 2 \abs{p}^2) \, (1 - 6 \abs{p}^2 + \abs{p}^4)
    - 4 (1 + 3\abs{p}^2) \langle w,p \rangle^2 \bigr) \, w\\[-1ex]
    & \qquad\qquad
      + \,2 (1-\abs{p}^2) \, (1+3\abs{p}^2) \langle w,p \rangle \, p
      \,-\, 4 (1 + 2\abs{p}^2) \, (1-\abs{p}^2) \, p\times w \Bigr).
\end{split}
\end{align}
Here, we use the standard cross product on $\R^3$ and the cross product on $\R^7$
that comes from the imaginary part of the product of two imaginary octonions.
It is straightforward but tedious to check that $x_0, y_0, x,y$ satisfy the
equations (\ref{realbrieskorn}).

On $\Sph^5\subset \R^3\times\R^3$ and on $\Sph^{13} \subset \R^7\times\R^7$
we consider the diagonal actions of $\SO(3)$ and $\Gtwo$, respectively.
The orbit spaces of these actions can again be identified with $D^2$ by the projection
maps $\Sph^{2n-1}\to D^2$,
$(p,w) \mapsto \abs{w}^2 - \abs{p}^2  - 2 i \,\langle p,w\rangle$.
Note that the preimage of the boundary of $D^2$ consists precisely of the
pairs $(p,w)$ for which $p$ and $w$ are linearly dependent.

\begin{thm}
The formulas (\ref{identbrieskorn}) above provide an $\SO(3)$-equivariant
diffeomorphism $\Sph^5 \to W^5_3$ and a $\Gtwo$-equivariant
diffeomorphism $\Sph^{13}\to W^{13}_3$.
\end{thm}

\begin{proof}
The maps $\psi:\Sph^{2n-1} \to W^{2n-1}_3$ defined by (\ref{identbrieskorn}) are
smooth, equivariant, and induce the identity between the orbit spaces $D^2$.
Hence they are homeomorphisms. Their inverses
$\psi^{-1}: W^{2n-1}_3 \to \Sph^{2n-1}$ can be computed explicitly:
The coefficients in the equations for $x$, $y$, and $x\times y$ as combinations
of $p$, $w$, and $p\times w$ are rational functions of $x_0$ and $y_0$
with nonzero denominators. The determinant of the coefficient matrix is a
polynomial of degree $12$ in $x_0$ and $y_0$ that can be seen to be
always greater than or equal to $ \tfrac{256}{9(3-2x_0)^8}$
if $x_0^2 + y_0^2 \le \tfrac{1}{4}$.
Hence, the coefficient matrix can be inverted even if $p$ and $w$ become
linearly dependent and $p$ and $w$ can be expressed as combinations
of $x$, $y$, and $x\times y$ with rational coefficients in $x_0$ and $y_0$
that do not have singularities for $x_0^2 + y_0^2 \le \tfrac{1}{4}$.
\end{proof}

If $n$ is different from $3$ and $7$ formula (\ref{identbrieskorn}) does not work.
What we need is to assign to $p,w\in \R^n$ a vector that is perpendicular to both
and that is different from zero if $p$ and $w$ are linearly independent.
Such a cross product exists only in dimensions $3$ and $7$ (see \cite{massey}).

In the rest of this section we will describe how formula (\ref{identbrieskorn})
was obtained in the case $n=3$. During this derivation we will meet a
simple formula for an injective map $\Sph^5\smallsetminus \{w = 0\} \to W^5_3$
that extends to a all odd dimensions and thus yields injective maps
$\Sph^{2n-1}\smallsetminus \{w = 0\} \to W^{2n-1}_3$.
These maps are given by substituting the expressions for $x$ and $y$ in
(\ref{identbrieskorn}) by
\begin{align*}
   -3x \,=\,& (\abs{p}^2 + 3\abs{w}^2 - 4 \langle \tfrac{w}{\abs{w}},p \rangle^2) p
        + 2 \abs{p}^2 \langle p, \tfrac{w}{\abs{w}} \rangle \tfrac{w}{\abs{w}}\, ,\\
    3y \,=\, & -(3\abs{p}^2 + \abs{w}^2) w + 6 \langle w,p \rangle p\,.
\end{align*}
The cross products in dimensions $3$ and $7$ are needed to twist these
maps such that they extend to diffeomorphisms $\Sph^5\to W^5_3$ and
$\Sph^{13} \to W^{13}_3$.

Now we start with the derivation of formula (\ref{identbrieskorn}).
There are two different para\-metrizations of $\Sigma^5\subset \Sigma^7$:
The explicit geodesic parametrization given in (\ref{lift}) describes a point
in $\Sigma^5$ by $(p,w) \in \Sph^5\subset \Imag\H \times \Imag\H$,
and the cohomogeneity one action of $\SO(3)\times \SO(2)$ on $\Sigma^5$
describes the same point by the parameter $s$ of the normal geodesic,
an angle $\theta$, and a unit quaternion $q\in\Sph^3$ (with several ambiguities).
More precisely, the identity
\begin{gather}
\label{paramident}
   q' \star \tilde \gamma_{\bmat p\\ w\emat}\bigl(\tfrac{\pi}{2}\bigr)
  \;=\; \bigl(\bmat \cos\theta & -\sin \theta\\ \sin\theta & \m \cos\theta\emat, q\bigr)
    \bullet \tilde \alpha(s)
\end{gather}
has to hold for some $q' \in \Sph^3$. In principle, one now solves for $s$, $\theta$,
and $\pm q$ in dependence of $p$ and $w$ (not caring about any ambiguities)
and plugs the results into the corresponding cohomogeneity one parametrization
\begin{align}
\begin{split}
\label{parambries}
    \bigl( \bmat \cos\theta & -\sin \theta\\ \sin\theta & \m \cos\theta\emat,
  & \pm q \bigr) \cdot \beta(s) =\\
  &  \bigl( -\tfrac{1}{2}\cos 2\theta \cos 2s - \tfrac{i}{2}\sin 2\theta \cos 2s,
    x(s,\theta,q) + i y(s,\theta,q) \bigr)
\end{split}
\end{align}
of $W^5_3$ where $x(s,\theta,q), y(s,\theta,q)\in \R^3$ will be evaluated below.
By Theorem\,\ref{equivdiff} it is clear that this procedure yields a well-defined smooth
diffeomorphism $\Sph^5 \to W^5_3$. The actual computations, however, are
lengthy and not straightforward. It thus seems appropriate to indicate how they
can be done efficiently. First, we identify $\R^3$ with the imaginary quaternions.
The homomorphism $\Sph^3 \to \SO(3)$ is then given by assigning
to $\pm q$ the matrix $(qi\bar q, qj\bar q, qk\bar q)\in \SO(\Imag\H)$.
With this identification $x(s,\theta,q), y(s,\theta,q)$ can be evaluated to
\begin{align}
\begin{split}
\label{xy}
    -3x(s,\theta,q) \,=\,
      & 2\,(1+\cos 2\theta\cos 2s)\, q(j\cos s\cos \theta - k \sin s\sin\theta)\bar q\\
      & \quad - (4 \cos 2\theta + \cos 2s)\, q(j\cos s\cos \theta + k \sin s\sin\theta)\bar q,\\
    3y(\theta,s) \,=\;
      & 2\,(1-\cos 2\theta\cos 2s)\, q(j\cos s\sin \theta + k \sin s\cos\theta)\bar q\\
      & \quad + (4 \cos 2\theta - \cos 2s)\, q(j\cos s\sin \theta - k \sin s\cos\theta)\bar q.
\end{split}
\end{align}
With a few computations one sees from (\ref{paramident}) that
\begin{gather*}
  \cos 2\theta = \tfrac{\abs{p}^2-\abs{w}^2}{
    \sqrt{(\abs{p}^2-\abs{w}^2)^2+4\langle p,w\rangle^2}}, \quad
  \sin 2\theta = \tfrac{2 \langle p,w\rangle}{
    \sqrt{(\abs{p}^2-\abs{w}^2)^2+4\langle p,w\rangle^2}},\\
  \cos 2s = \sqrt{(\abs{p}^2-\abs{w}^2)^2+4\langle p,w\rangle^2}, \quad
  \sin 2s = 2 \sqrt{(\abs{p}^2 \abs{w}^2-\langle p,w\rangle^2}.
\end{gather*}
Moreover, (\ref{paramident}) is equivalent to
\begin{gather*}
   \bmat \x\m\cos\theta & \m \sin \theta\\ \x-\sin\theta & \m \cos\theta\emat
   \cdot \tilde \gamma_{\bmat p\\ w\emat}\bigl(\tfrac{\pi}{2}\bigr)
  \;=\; \bar q' \star (q \bullet \tilde \alpha(s)).
\end{gather*}
Computing formally the ``determinant'' $ad-bc$ of the quaternionic $2\times 2$
matrices $\bmat a & c\\ b & d\emat$ on both sides of the latter equation one obtains
\begin{gather*}
  \bar q' \cos 2s =  \Bigl( (\abs{p}^2-\abs{w}^2)\tfrac{w}{\abs{w}}
    - 2\langle \tfrac{w}{\abs{w}},p \rangle p \Bigr) e^{\frac{\pi}{2}p} \tfrac{\bar w}{\abs{w}} \,q.
\end{gather*}
This identity can now be plugged into (\ref{paramident}) and the result allows us
to evaluate (\ref{xy}):
\begin{align*}
   -3x \,=\,& \tfrac{w}{\abs{w}} e^{-\frac{\pi}{2}p}
      \Bigl( (\abs{p}^2 + 3\abs{w}^2 - 4 \langle \tfrac{w}{\abs{w}},p \rangle^2) p
        + 2 \abs{p}^2 \langle p, \tfrac{w}{\abs{w}} \rangle \tfrac{w}{\abs{w}} \Bigr)
      e^{\frac{\pi}{2}p} \tfrac{\bar w}{\abs{w}}\, ,\\
    3y \,=\, & \tfrac{w}{\abs{w}} e^{-\frac{\pi}{2}p}
    \Bigl( -(3\abs{p}^2 + \abs{w}^2) w + 6 \langle w,p \rangle p \Bigr)
      e^{\frac{\pi}{2}p} \tfrac{\bar w}{\abs{w}}\,.
\end{align*}
Expressing all quaternionic products by inner products and cross products
we obtain the formulas
\begin{align*}
  3 x \,=\,& \bigl( 3 - 2 \abs{p}^2
      - 2 \, \tfrac{1+\cos \pi\abs{p}}{1-\abs{p}^2} \langle w,p \rangle^2 \bigr) \, p\\
    & \quad - 2 \bigl( 3 + \abs{p}^2 \tfrac{1+\cos \pi\abs{p}}{1-\abs{p}^2}
      - 2 \, \tfrac{1+\cos \pi\abs{p}}{(1-\abs{p}^2)^2} \langle w,p \rangle^2 \bigr)
      \langle p,w \rangle \, w\\
    &\qquad - 2 \abs{p}^2 \tfrac{\sin \pi\abs{p}}{(1-\abs{p}^2)\abs{p}}
    \, \langle w,p \rangle \, p\times w,\\
  3y \,=\,& - ( 1+ 2 \abs{p}^2) \cos \pi\abs{p} \cdot w 
      \,+\, 2\, \tfrac{-1 + 4 \abs{p}^2
      + (1+2\abs{p}^2) \cos \pi\abs{p}}{\abs{p}^2(1-\abs{p}^2)}
      \langle w,p \rangle^2 w\\
    & \quad - \tfrac{-1 + 4 \abs{p}^2
      + (1+2\abs{p}^2) \cos \pi\abs{p}}{\abs{p}^2(1-\abs{p}^2)}
      \, \abs{w}^2 \langle w,p \rangle \, p
      \,-\, (1 + 2 \abs{p}^2) \tfrac{\sin \pi\abs{p}}{\abs{p}} \, p\times w
 \end{align*}
where all the fractions are real analytic functions of $\abs{p}$.
This can now be seen as a final formula for the diffeomorphism $\Sph^5\to W^5_3$.
In formula (\ref{identbrieskorn}) we passed to an isotopic rational version by substituting
$\sin \tfrac{\pi}{2}\abs{p}$ and $\cos \tfrac{\pi}{2}\abs{p}$ by $\tfrac{2\abs{p}}{1+\abs{p}^2}$ and
$\tfrac{1-\abs{p}^2}{1+\abs{p}^2}$, respectively.

\begin{rem}
Note that the diffeomorphisms of Theorem\,\ref{nonlinear} equip $W^5_3$
and $W^{13}_3$ with explicit $\SO(3)$ and $\Gtwo$ invariant metrics of
constant curvature $1$. Wilking (unpublished) proved that there do not exist
$\SO(n)$-invariant metrics with positive sectional curvature on any
of the $W^{2n-1}_d$ with $n > 3$ and odd $d > 1$. Moreover, it was
shown \cite{gvwz} that there do not exist cohomogeneity one metrics with
nonnegative sectional curvature on any of the $W^{2n-1}_d$ with $n > 3$
and odd $d > 1$.
\end{rem}

\bigskip

\section{Nonlinear cohomogeneity one actions on Euclidean spheres}
\label{nlEucl}
In this section we present the first explicit formulas for cohomogeneity one actions
of $\OO(2) \times \SO(3)$ and $\OO(2)\times \Gtwo$ on the Euclidean spheres
$\Sph^5$ and $\Sph^{13}$ that are equivalent to the standard cohomogeneity one
actions  on the Brieskorn spheres $W^5_3$ and $W^{13}_3$ (see section\,\ref{briesid}).

The essential parts of these actions are the nonlinear subactions of $\SO(2) \subset \OO(2)$.
It is convenient, however, to describe the linear parts first:
Let $p$ and $w$ denote two imaginary quaternions (octonions)
with $\abs{p}^2 + \abs{w}^2 = 1$. The action of $\SO(3)= \Sph^3/\{\pm 1\}$
on $\Sph^5$ is given by
\begin{gather*}
  \SO(3) \times \Sph^5 \to \Sph^5, \quad
  (\pm q) \bullet \bsmat p\\ w\esmat = q \bsmat p \\ w \esmat \bar q
  = \bsmat q p \bar q\\ q w \bar q\esmat.
\end{gather*}
The $\Gtwo$-action on $\Sph^{13}$ is defined in the same diagonal way.
(Recall that $\Gtwo$ is the automorphism group of the octonions).
The $\OO(2)$-actions on $\Sph^5$ and $\Sph^{13}$ contain the
linear $\Z_2$-subactions
\begin{gather*}
  \Z_2 \times \Sph^5 \to \Sph^5, \quad
  \bsmat 1 & \m 0\\ 0 & -1 \esmat \bullet \bsmat p\\ w\esmat =
  \bsmat \;\, p\\ \x -w\esmat
\end{gather*}

We will now turn to the nonlinear $\SO(2)$-actions. In order to write them
down explicitly we need some preparatory work.
Let $e^p$ denote the exponential map of $\Sph^3\subset \H$ (or $\Sph^7 \subset \Ca$).
For $\theta \in \R$ set
\begin{gather*}
  \bsmat p_{\theta} \\ w_{\theta} \esmat :=
  \bsmat \cos \theta & -\sin \theta\\ \sin \theta & \m \cos \theta \esmat \cdot
  \bsmat p\\ w\esmat =
  \bsmat p \cos \theta - w \sin \theta\\ p \sin \theta + w\cos \theta\esmat
\end{gather*}
and
\begin{align*}
  Q\bigl(\bsmat p\\ w\esmat, \theta \bigr) = 
  \tfrac{w}{\abs{w}} e^{-\frac{\pi}{2} p}
    \tfrac{w_{\theta}}{\abs{w_{\theta} }} \tfrac{\bar w}{\abs{w}}
    e^{\frac{\pi}{2}p_{\theta}}
    \tfrac{\bar w_{\theta} }{\abs{w_{\theta}}}.
\end{align*}
At first glance one would not expect that this formula defines
a smooth map.

\begin{lem}
$Q$ extends to an analytic map $\Sph^5\times \Sph^1 \to \Sph^3$ and
$\Sph^{13} \times \Sph^1 \to \Sph^7$, respectively.
\end{lem}

\begin{proof}
Expanding the exponential maps in the definition of $Q$ and applying the
two identities $p_{\theta} = w_{(\theta+\frac{\pi}{2})}$ and
$w_{\theta} w_{\kappa} w_{\tau} = w_{\tau} w_{\kappa} w_{\theta}$ one obtains
\begin{align*}
  Q\bigl(\bsmat p\\ w\esmat, \theta \bigr) =
  & \; w w_{\theta} \bar w \bar w_{\theta} \cdot \tfrac{\cos \frac{\pi}{2}\abs{p}}{1-\abs{p}^2}
    \cdot \tfrac{\cos \frac{\pi}{2}\abs{p_{\theta}}}{1-\abs{p_{\theta}}^2}
  \,-\, p_{\theta} p \cdot \tfrac{\sin \frac{\pi}{2}\abs{p}}{\abs{p}}
  \cdot \tfrac{\sin \frac{\pi}{2}\abs{p_{\theta}}}{\abs{p_{\theta}}} \\
  & + w p_{\theta} \bar w \cdot \tfrac{\cos \frac{\pi}{2}\abs{p}}{1-\abs{p}^2} \cdot
    \tfrac{\sin \frac{\pi}{2}\abs{p_{\theta}}}{\abs{p_{\theta}}}
  \,-\, w_{\theta} p \bar w_{\theta} \cdot \tfrac{\sin \frac{\pi}{2}\abs{p}}{\abs{p}} \cdot
    \tfrac{\cos \frac{\pi}{2}\abs{p_{\theta}}}{1-\abs{p_{\theta}}^2} \qedhere
\end{align*}
\end{proof}

\begin{lem}
\label{actprop}
$Q$ has the following property:
\begin{gather*}
  Q\bigl(\bsmat p\\ w\esmat, \theta \bigr)\,
    Q\bigl(\bsmat p_{\theta}\\ w_{\theta}\esmat, \tau \bigr)
  = Q\bigl(\bsmat p\\ w\esmat, \theta+\tau \bigr).
\end{gather*}
\end{lem}

\begin{proof}
This property is based on the identity 
$w_{\theta} \bar w w_{\tau} = w_{\tau} \bar w w_{\theta}$.
\end{proof}

\begin{thm}
\label{nonlinear}
The assignment
\begin{gather*}
  \bsmat \cos \theta & -\sin \theta\\ \sin \theta & \m \cos \theta \esmat \bullet
  \bsmat p\\ w\esmat :=
  Q \bigl(\bsmat p\\ w\esmat, \theta \bigr)\,
    \bsmat p_{\theta}\\ w_{\theta}\esmat\,
  \overline{ Q\bigl(\bsmat p\\ w\esmat, \theta \bigr)}.
\end{gather*}
defines nonlinear $\SO(2)$-actions on $\Sph^5$ and $\Sph^{13}$
that extend to cohomogeneity one actions of $\OO(2)\times \SO(3)$
and $\OO(2) \times \Gtwo$, respectively. These latter actions are equivalent to
the standard actions on the Brieskorn spheres $W^5_3$ and~$W^{13}_3$.
\end{thm}

\begin{proof}
The map $Q: \Sph^5\times \Sph^1 \to \Sph^3$ is equivariant under conjugation
with unit quaternions, i.e.,
\begin{gather*}
  Q\bigl(\bsmat qp\bar q\\ qw\bar q\esmat, \theta \bigr)
  \,=\, q \, Q\bigl(\bsmat p\\ w\esmat, \theta \bigr) \,\bar q,
\end{gather*}
and $Q:\Sph^{13} \times \Sph^1 \to \Sph^7$ is in an analogous way
equivariant under $\Gtwo$.
With Lemma\,\ref{actprop} it is now straightforward to check that the assignment of
Theorem\,\ref{nonlinear} defines $\SO(2)$-actions on $\Sph^5$ and $\Sph^{13}$
which commute with the $\SO(3)$-action on $\Sph^5$ and the $\Gtwo$-action
on $\Sph^{13}$. It can be proved in various ways that $\Sph^5$ and $\Sph^{13}$
equipped with the full actions of $\OO(2)\times \SO(3)$ and $\OO(2) \times \Gtwo$
are equivariantly diffeomorphic to $W^5_3$ and~$W^{13}_3$, e.g., by computing
the isotropy groups along the curve
\begin{gather*}
  s \mapsto \bsmat j \cos s\\ ( k \cos(\pi\cos s) - i \sin(\pi\cos s)) \sin s\esmat
\end{gather*}
which corresponds precisely to the geodesics $\alpha$ on $\Sigma^5$ and
$\beta$ on $W^5_3$ and $W^{13}_3$ under the identifications established
in the previous sections. (For the isotropy group computation note that if
$p$ and $w$ anticommute and have the same norm then we have
$p_{\theta} = e^{-v\frac{\theta}{2}} p e^{v\frac{\theta}{2}}$
where $v = \frac{p}{\abs{p}}\frac{w}{\abs{w}}$
and a similar expression for $w_{\theta}$.)
\end{proof}

\begin{rem}
The formula of Theorem\,\ref{nonlinear} was obtained by pulling back the
$\bullet$-action on $\Sigma^5$ by the explicit diffeomorphism $\Sph^5\to\Sigma^5$
given in (\ref{expldiffeo}).
\end{rem}

\begin{rem}
For $\theta = \pi$ the formula of Theorem\,\ref{nonlinear} gives exotic involutions
on $\Sph^5$ and $\Sph^{13}$. These are studied in the paper \cite{involutions}.
\end{rem}

\begin{rem}
If one substitutes $\tfrac{\pi}{2}$ in the definition of $Q$ by $(2m+1)\tfrac{\pi}{2}$
then one obtains an action that is conjugate to the original action by $\sigma^m$
where $\sigma$ is the restriction of the exotic diffeomorphism $\sigma:\Sph^6\to \Sph^6$
to $\Sph^5$ (see \cite{duran}).
\end{rem}

\bigskip

%
%

\end{document}